\documentclass[12pt]{article}

\usepackage{amssymb}
\usepackage{url}
\usepackage{gastex}

\newenvironment{Proof}{\rm \trivlist \item[\hskip \labelsep{\bf
Proof.}]}{\cqfd\endtrivlist}

\def\cqfd{\skip10=\parfillskip\parfillskip=0pt
\enspace\hfill\symbolecqfd\par\parfillskip=\skip10\par\medskip}
\def\symbolecqfd{\rlap{$\sqcap$}$\sqcup$}
\def\preuve{\begin{Proof}}
\def\proof{\begin{Proof}}
\def\eop{\end{Proof}}

\newtheorem{theorem}{Theorem}[section]
\newtheorem{proposition}[theorem]{Proposition}
\newtheorem{lemma}[theorem]{Lemma}
\newtheorem{corollary}[theorem]{Corollary}

\newtheorem{pro-example}[theorem]{Example}
\newenvironment{example}{\begin{pro-example}\rm}{\cqfd\end{pro-example}}

\newtheorem{pro-remark}[theorem]{Remark}
\newenvironment{remark}{\begin{pro-remark}\rm}{\cqfd\end{pro-remark}}

\newtheorem{pro-fact}[theorem]{Fact}
\newenvironment{fact}{\begin{pro-fact}\rm}{\cqfd\end{pro-fact}}

\newcommand{\calP}{\mathcal{P}}

\let\phi\varphi

\def\lefg{\le_{\textsf{\scriptsize fg}}}

\def\capac{\textsf{cap}}
\def\card{\textsf{card}}
\def\deg{\textsf{deg}}
\def\rank{\textsf{rank}}
\def\Aut{\textsf{Aut}}
\def\hl{\textsf{hl}}
\def\link{\textsf{link}}
\def\red{\textsf{red}}

\def\W{\mathbb{W}}
\def\Z{\mathbb{Z}}

\def\cc{\textsf{cc}}
\def\id{\textsf{id}}

\def\inv{^{-1}}

\begin{document}

\title{On the complexity of the Whitehead minimization
problem\protect\footnote{Part of this work was produced while the
second author was a Visiting Professor at LaBRI, Universit\'e
Bordeaux-1.  The second author gratefully acknowledges partial support
by DGI (Spain) through grant BFM2003-06613.  The third author
gratefully acknowledges support from the ESF programme
\textsc{AutoMathA}.}}

\markright{\protect Complexity of the Whitehead minimization problem}

\author{Abd\'o Roig
\small{(\url{abdo.roig@upc.edu})}
% \protect\footnote{Universitat Polit\`ecnica de Catalunya}
\\
\small{Universitat Polit\`ecnica de Catalunya}\\
\\ Enric Ventura
\small{(\url{enric.ventura@upc.edu})}\protect\footnote{EPSEM,
Universitat Polit\`ecnica de Catalunya -- Av.  Bases de Manresa 61 --
73 08242 Manresa, Barcelona -- Spain} \\
\small{Universitat Polit\`ecnica de Catalunya}\\
\\ Pascal Weil
\small{(\url{pascal.weil@labri.fr})}\protect\footnote{LaBRI -- 351
cours de la Lib\'eration -- 33405 Talence Cedex -- France} \\
\small{LaBRI (CNRS and Universit\'e Bordeaux-1)}}

\date{}

\maketitle

\begin{abstract}
    The Whitehead minimization problem consists in finding a minimum
    size element in the automorphic orbit of a word, a cyclic word or
    a finitely generated subgroup in a finite rank free group.  We
    give the first fully polynomial algorithm to solve this problem,
    that is, an algorithm that is polynomial both in the length of the
    input word and in the rank of the free group.  Earlier algorithms
    had an exponential dependency in the rank of the free group.  It
    follows that the primitivity problem -- to decide whether a word
    is an element of some basis of the free group -- and the free
    factor problem can also be solved in polynomial time.
\end{abstract}

\bigskip

%%%%%%%%%%%%%%%%%%%%%%%
%%%%%%%%%%%%%%%%%%%%%%%

Let $F$ be a finite rank free group and let $A$ be a fixed (finite)
basis of $F$.  The elements of $F$ can be represented by reduced words
over the symmetrized alphabet $A\cup\bar A$, and the finitely
generated subgroups of $F$ by certain finite graphs whose edges are
labeled by letters from $A$ (obtained by the technique of so-called
Stallings foldings \cite{Stallings}, see \cite{KM} and
Section~\ref{sec 1}).  We measure the size of an element of $F$ by the
length of the reduced word representing it, and the size of a finitely
generated subgroup of $F$ by the number of vertices of the graph
representing it.  The \textit{Whitehead minimization problem} consists
in finding a minimum size element in the automorphic orbit of a given
word or a given finitely generated subgroup.  An important variant
considers rather as input conjugacy classes of words (the so-called
cyclic words) or subgroups.

As we will see (Section~\ref{sec W theorem}), the minimization problem
for words, cyclic words and subgroups can be reduced to the problem
for conjugacy classes of finitely generated subgroups, so we will
often discuss only the latter problem.

The Whitehead minimization problem is a fragment of the larger
\textit{equivalence problem}, where one must decide, given two
subgroups of $F$, whether they sit in the same automorphic orbit.
More precisely, in view of a result of Gersten \cite{Gersten} (which
generalizes to subgroups a classical result of Whitehead that applies
to words \cite{Whitehead}, see \cite[Sec.  I.4]{LS}),
the first part of the (only known) algorithm to solve the equivalence
problem consists in finding minimum size elements in the automorphic
orbits of the given subgroups, that is, in solving the Whitehead
minimization problem for these two subgroups.

Moreover, any solution of the Whitehead minimization problem for words
implies a solution of the \textit{primitivity problem}: to decide
whether a given word is an element of some basis of $F$.  Indeed, a
word is primitive if and only if the minimum length of an element in
its automorphic orbit is 1.  Similarly, a subgroup is a free factor of
$F$ if and only if the minimum size of an element in its automorphic
orbit is 1, so any solution of the Whitehead minimization problem for
subgroups implies a solution of the \textit{free factor problem}: to
decide whether a given subgroup is a free factor of $F$.

As hinted above, a classical algorithm is known to solve the Whitehead
minimization problems.  The algorithm for the word case is based on
Whitehead's theorem \cite{Whitehead} (see \cite[Prop.  I.4.17]{LS}),
and the algorithm for the subgroup case relies on a generalization of
this theorem due to Gersten \cite[Corol.  2]{Gersten} (Theorem~\ref{W
theorem} below).  Both algorithms are polynomial in $n$, the size of
the input, and exponential in $r$, the rank of the free group $F$, see
Section~\ref{sec W automorphisms} below and \cite{ADAG} for a more
detailed analysis.

In a recent paper, Haralick, Miasnikov and Myasnikov \cite{HMM} (see
also Miasnikov and Myasnikov \cite{ADAG}) present a number of
heuristics and experiments on different implementation strategies for
the algorithm regarding words, that tend to show that the actual
dependency of the problem in the rank of $F$ is much lower than
exponential, at least in the generic case.  Haralick and Miasnikov
\cite{ADMH} pursue in that direction by giving a polynomial-time
stochastic algorithm for the same problem.  Another recent paper, by
Silva and Weil \cite{SilvaWeil} gives an exact algorithm for solving
the free factor problem on input $H$, which is polynomial in the size
of $H$ and exponential in $\rank(F) - \rank(H)$.

The main result of this paper confirms the intuition of
\cite{HMM,ADAG,ADMH} and improves the complexity result in
\cite{SilvaWeil}.  Indeed, we give a fully polynomial solution to the
Whitehead minimization problem, that is, an exact algorithm that is
polynomial in both the size of the input and the rank of $F$.
Interestingly, this result is obtained with only a small amount of new
mathematical results.  Our algorithm is in fact a minor modification
of the classical Whitehead method (the algorithm mentioned above), and the study
of its complexity relies on the conjunction of three ingredients:
\begin{itemize}
    \item[(1)] a representation of a Whitehead automorphism as a
    pointed cut of the set $A\cup\bar A$, that is, a partition of
    $A\cup\bar A$ into two disjoint subsets $Y$ and $Z$ and the choice
    of a letter $a\in Y$ such that $\bar a\in Z$,
    
    \item[(2)] an exact computation of the effect of such an
    automorphism on the length of a cyclic word $u$ (resp.  the size
    of a conjugacy class $H$ of finitely generated subgroups) in terms
    of the capacity of the associated cut on the Whitehead graph
    associated with $u$ (resp.  a generalization of the Whitehead
    graph which we call the Whitehead hypergraph of $H$),
    
    \item[(3)] and an algorithmic complexity result on finding a
    minimum capacity cut in a graph (resp.  a hypergraph).
\end{itemize}

The first ingredient is classical in combinatorial group theory (see
for instance \cite[Prop.  I.4.16]{LS}).  The second ingredient can be
described as a rewording of a formula of Gersten \cite[Corol.
2]{Gersten} proved in \cite[Prop.  10.3]{Kalajd}.  And the last one
can be reduced to standard problems in combinatorial optimization, for
which there exist several polynomial-time algorithms in the
literature.

In Section~\ref{sec 1}, we fix the notation to handle words, cyclic
words and subgroups of $F$, and to describe Whitehead automorphisms.
We also discuss the foundational theorem of Whitehead and its
generalization by Gersten, that gives rise to the known algorithm
solving the Whitehead minimization problem.  As indicated above this
algorithm is polynomial in the size of the input and exponential in
the rank of the ambient free group.

In Section~\ref{sec choice}, we introduce the Whitehead hypergraph of
a cyclically reduced subgroup, and we present a rewording of a formula
of Gersten, describing the effect of a Whitehead automorphism on the
size of a conjugacy class of finitely generated subgroups (and
generalizing a classical result of Whitehead on cyclic words, see
\cite[Prop.  I.4.16]{LS}).  This technical analysis helps show how the
exponential dependency in the usual Whitehead algorithm can be
removed, provided a certain graph-theoretic problem, namely the
min-cut problem for undirected hypergraphs, can be solved in
polynomial time.

We discuss existing polynomial-time algorithms to solve the min-cut
problem in Section~\ref{sec main}, thus completing our proof.
Finally, we consider some consequences of our main result.

We conclude this introduction with a remark on complexity computation.
Since the rank $r$ of $F$ is part of the input, we consider complexity
functions under the \textit{bit cost assumption}: $r$ is the
cardinality of $A$ and each letter is identified by a bit string of
length $\log r$, so that the elementary operations on $A$ (reading or
writing a letter, comparing two letters, etc) require $O(\log r)$
units of time.

%%%%%%%%%%%%%%%%%%%%%%
\section{The Whitehead minimization problem}\label{sec 1}

In this paper, $F$ denotes a finitely generated free group and $A$
denotes a fixed basis of $F$.  We let $r = \rank(F) = \card(A)$.

%%%%%%%%%%%%%%%%%%%%%%
\subsection{Words, graphs and subgroups}

Elements of $F$ can be represented as usual by \textit{reduced words}
on the symmetrized alphabet $\tilde A = A \cup \bar A$, and we write
$u\in F(A)$ to indicate that the element $u$ of $F$ is given by a
reduced word on alphabet $\tilde A$.

Recall that the mapping $a \mapsto \bar a$ is extended to the set of
all words over the alphabet $\tilde A$ by letting $\bar{\bar a} = a$
for each $a\in A$, and $\overline{ua} = \bar a\,\bar u$ for each word
$u$ and each letter $a\in\tilde A$.

%%%%%%%%%%%%%%%%%%%%%%
\subsubsection{Graphs}

To represent finitely generated subgroups of $F$, we use finite
$A$-graphs.  An \textit{$A$-graph} is a directed graph, whose edges
are labeled by letters in $\tilde A$.  More precisely, it is a pair
$\Gamma = (V,E)$ with $E \subseteq V \times \tilde A \times V$.  The
elements of $V$ are called \textit{vertices}, the elements of $E$ are
called \textit{edges}, and we say that there is an edge from $x$ to
$y$ labeled $a$ if $(x,a,y)\in E$.  We denote respectively by
$\alpha$, $\lambda$ and $\omega$ the first, second and third
coordinate projections from $E$ to $V$, $\tilde A$ and $V$.  The map
$\lambda$ is called the \textit{labeling function}.

We measure the \textit{size} of an $A$-graph $\Gamma$ by the number of
its vertices, and we write $|\Gamma| = \card(V)$.

A \textit{dual} $A$-graph is one in which for each $a\in A$, there is
an edge from vertex $x$ to vertex $y$ labeled $a$ if and only if there
is one from $y$ to $x$ labeled $\bar a$.  That is, $(x,a,y)\in E$ if
and only if $(y,\bar a,x)\in E$.

A dual $A$-graph is \textit{reduced} if whenever there are $a$-labeled
edges from $x$ to $y$ and from $x'$ to $y$, then $x = x'$ ($a\in
\tilde A$, $x,x',y \in V$).  It is \textit{cyclically reduced} if it
is reduced and, for each vertex $x$, there exist at least 2 edges into
$x$.

If $\Gamma$ is an $A$-graph and $x\in V$, we define the \textit{link}
of $x$ to be the set of edges into $x$, and the \textit{hyperlink} of 
$x$ to be the set of labels of these edges,
$$\link_\Gamma(x) = \{e\in E \mid \omega(e) = x\},\quad \hl_\Gamma(x) = 
\{\lambda(e) \mid e\in\link_\Gamma(x)\}.$$
(The reason for the terminology \textit{hyperlink} will become clear 
in Section~\ref{sec hypergraph}.)

Let $\Gamma$ be a dual $A$-graph.  By an immediate rewording of the
definitions, we see that $\Gamma$ is reduced if and only if $\lambda$
establishes a bijection from $\link_\Gamma(x)$ to $\hl_\Gamma(x)$; and
$\Gamma$ is cyclically reduced if in addition, $\card(\hl_\Gamma(x))
\ge 2$ for each $x\in V$.  A vertex $x$ such that
$\card(\hl_\Gamma(x)) \le 1$, is called an \textit{endpoint} of
$\Gamma$.

We say that an $A$-graph is \textit{connected} if the underlying
undirected graph is connected.  In the case of a dual $A$-graph, this
is the case if and only if, for any vertices $x,y$, there exists a
path from $x$ to $y$.

The following particular case will be important for our purpose.  A
\textit{circular} $A$-graph is a connected dual $A$-graph in which the
link of each vertex has exactly 2 elements.  A cyclically reduced
circular $A$-graph is called a \textit{cyclic word}.

\begin{example}\label{example graphs}
    When representing $A$-graphs, we draw only the positively labeled
    edges, that is, those labeled by letters of $A$.  Figure~\ref{fig
    graphs} shows such $A$-graphs.
    In that figure, $\Gamma_1$, $\Gamma_2$ and $\Gamma_3$ are reduced,
    $\Gamma_0$ is not.  Vertices 2 and 3 are endpoints in $\Gamma_2$
    and $\Gamma_3$ respectively, while $\Gamma_1$ has no endpoint.  In
    $\Gamma_2$, the hyperlink of vertex 1 is $\{a, \bar a, \bar b\}$
    and the hyperlink of vertex 2 is $\{\bar a\}$.
\end{example}    

\begin{figure}[ht]
    \centering
    \begin{picture}(126,65)(0,-65)

    \node[Nw=5.0,Nh=5.0,Nmr=1.0](n2)(0.0,-8.0){}

    \node[Nw=5.0,Nh=5.0,Nmr=1.0](n3)(22.0,-8.0){1}

    \node[Nw=5.0,Nh=5.0,Nmr=1.0](n4)(44.0,-8.0){}

    \node[Nw=5.0,Nh=5.0,Nmr=1.0](n5)(11.0,-24.0){}

    \node[Nw=5.0,Nh=5.0,Nmr=1.0](n6)(33.0,-24.0){}
    
    \put(48.0,-24.0){$\Gamma_0$}

    \node[Nw=5.0,Nh=5.0,Nmr=1.0](n22)(77.0,-8.0){1}

    \node[Nw=5.0,Nh=5.0,Nmr=1.0](n23)(99.0,-8.0){}

    \node[Nw=5.0,Nh=5.0,Nmr=1.0](n24)(88.0,-24.0){}

    \put(103.0,-24.0){$\Gamma_1$}

    \node[Nw=5.0,Nh=5.0,Nmr=1.0](n25)(22.0,-40.0){1}

    \node[Nw=5.0,Nh=5.0,Nmr=1.0](n26)(44.0,-40.0){}

    \node[Nw=5.0,Nh=5.0,Nmr=1.0](n27)(33.0,-56.0){}

    \node[Nw=5.0,Nh=5.0,Nmr=1.0](n34)(0.0,-40.0){2}

    \put(48.0,-56.0){$\Gamma_2$}

    \node[Nw=5.0,Nh=5.0,Nmr=1.0](n28)(77.0,-40.0){}

    \node[Nw=5.0,Nh=5.0,Nmr=1.0](n29)(99.0,-40.0){}

    \node[Nw=5.0,Nh=5.0,Nmr=1.0](n30)(88.0,-56.0){}

    \node[Nw=5.0,Nh=5.0,Nmr=1.0](n35)(121.0,-40.0){3}

    \put(103.0,-56.0){$\Gamma_3$}

    \drawedge[ELside=r](n3,n2){$a$}

    \drawedge(n3,n4){$a$}

    \drawedge(n3,n5){$b$}

    \drawedge[ELside=r](n3,n6){$b$}

    \drawedge(n5,n2){$b$}

    \drawedge(n4,n6){$a$}

    \drawedge(n22,n23){$a$}

    \drawedge[ELside=r](n22,n24){$b$}

    \drawedge[curvedepth=3.0](n23,n24){$a$}

    \drawedge[curvedepth=3.0](n24,n23){$b$}

    \drawedge(n34,n25){$a$}

    \drawedge(n25,n26){$a$}

    \drawedge[curvedepth=3.0](n26,n27){$a$}

    \drawedge[curvedepth=3.0](n27,n26){$b$}

    \drawedge[ELside=r](n25,n27){$b$}

    \drawedge(n28,n29){$a$}

    \drawedge(n29,n35){$b$}

    \drawedge[ELside=r](n28,n30){$b$}

    \drawedge[curvedepth=3.0](n30,n29){$b$}

    \drawedge[curvedepth=3.0](n29,n30){$a$}

\end{picture}
\caption{Four $A$-graphs}
\label{fig graphs}
\end{figure}
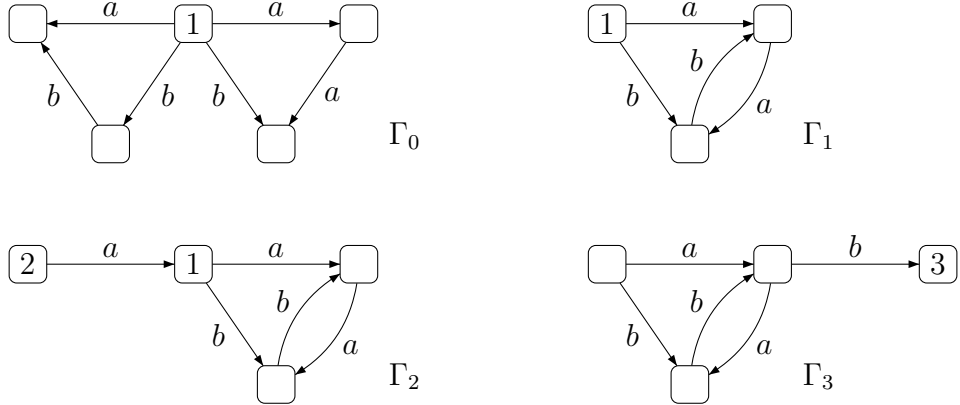

%%%%%%%%%%%%%%%%%%%%%%%
\subsubsection{Reduction}

Let $\Gamma = (V,E)$ be a dual $A$-graph, and let $x,y$ be distinct
vertices.  The $A$-graph obtained from $\Gamma$ \textit{by identifying
vertex $y$ to vertex $x$} is constructed as follows: its vertex set is
$V\setminus\{y\}$; and its edge set is obtained from $E$ by replacing
everywhere $y$ by $x$.  The resulting $A$-graph is again dual. Note 
that identifying $y$ to $x$ or $x$ to $y$ yields isomorphic $A$-graphs.

If $\Gamma$ is not reduced, there exist pairs of distinct edges
$(x,a,z)$ and $(y,a,z)$ (that is, edges with the same label that point
to the same vertex).  An \textit{elementary reduction} of $\Gamma$ is
the $A$-graph that results from identifying vertex $y$ to vertex $x$
in such a situation (which automatically implies the identification of
the two distinct edges).

The \textit{reduction of a dual $A$-graph} $\Gamma$ consists in
repeatedly performing elementary reductions, as long as some are
possible.  This is a terminating process since we consider only finite
graphs, and each elementary reduction properly decreases the size of
the graph.  The resulting graph, denoted by $\red(\Gamma)$, is
reduced, and it is well known that it does not depend on the sequence
of elementary reductions that were performed (that is: the process of
elementary reductions is confluent, see \cite{Stallings}).

If $\Gamma$ is reduced, an \textit{elementary trimming} consists in
removing an endpoint and the edges adjacent to it.  The \textit{cyclic
reduction} of $\Gamma$ consists in repeatedly applying elementary
trimmings, as long as it is possible.  The resulting graph is
cyclically reduced, it does not depend on the sequence of elementary
trimmings performed, and it is called the \textit{cyclic core} of
$\Gamma$, written $\cc(\Gamma)$.  Clearly, $\cc(\Gamma)$ is a subgraph
of $\Gamma$.

If $x$ is a vertex of $\Gamma$, there exists a unique shortest path
from $x$ to a vertex of $\cc(\Gamma)$.  We call this path the
\textit{branch} of $\Gamma$ at $x$, and we denote by $b(x)$ the label
of that path, and by $\beta(x)$ its extremity in $\cc(\Gamma)$.  If
$x$ is already in $\cc(\Gamma)$, then $\beta(x) = x$ and $b(x)$ is the
empty word, $b(x) = 1$.

\begin{example}\label{example reduction}
    With reference to the $A$-graphs in Example~\ref{example graphs} 
    and Figure~\ref{fig graphs}, we have
    $$\red(\Gamma_0) = \Gamma_1 = \cc(\Gamma_2) = \cc(\Gamma_3).$$
    In graph $\Gamma_2$, we have $\beta(2) = 1$ and $b(2) = a$.
\end{example}

Let $u = a_1\cdots a_n\in F(A)$.  The word $u$ is said to be
\textit{cyclically reduced} if $a_n \ne \bar a_1$, if and only if the
word $u^2$ is reduced.  It is a standard observation that, in general,
there exists a unique cyclically reduced word $w$ such that $u =
vw\bar v$ for some $v\in F(A)$.  The word $w$ is called the
\textit{cyclic core} of $u$, written $\cc(u)$.  Let $\Gamma(u)$ be the
circular graph with vertex set $\Z/n\Z = \{1,\ldots,n\}$ and with
edges $(i,a_i,i+1)$ and $(i+1,\bar a_i,i)$ for each $1\le i\le n$.
Reducing $\Gamma(u)$ yields the graph shown in Figure~\ref{fig
handle}, and $\cc(\Gamma(u)) = \Gamma(\cc(u))$.

\begin{figure}[ht]
    \centering
    \begin{picture}(47,22)(0,-24)

    \node[linewidth=0.0,Nw=6.0,Nh=6.0,Nmr=1.0](n0)(0.0,-15.0){1}

    \node[linewidth=0.0,Nw=6.5,Nh=6.5,Nmr=1.0](n1)(30.0,-15.0){$\scriptstyle
    \beta(1)$}

    \drawedge(n0,n1){$v = b(1)$}

    \drawloop[loopdiam=22.0,loopangle=0.0](n1){$w = \cc(u)$}
    
    \end{picture}
    \caption{The graph $\red(\Gamma(u))$, where $u = vw\bar v$, $v =
    b(1)$ and $w = \cc(u)$}
    \label{fig handle}
\end{figure}
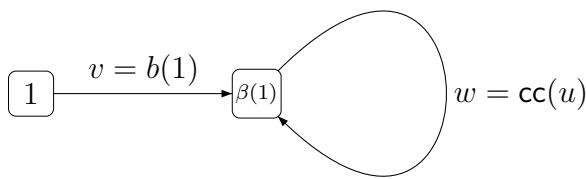

\begin{fact}\label{fact reducing and trimming}
    It is verified in \cite[Sec.  3.2]{SilvaWeil} that reducing and
    trimming an $n$-vertex dual $A$-graph takes time $O(n^2\log(nr))$.
    See Touikan \cite{Touikan} for a faster algorithm.
    
    If $u\in F(A)$ has length $n$, its cyclic core $\cc(u)$ is
    computed in time $O(n\log r)$.  In particular, reducing and
    trimming an $n$-vertex circular graph takes time $O(n\log r)$.
\end{fact}    

%%%%%%%%%%%%%%%%%%%%%%%
\subsubsection{Graphs and subgroups}\label{sec subgroups}

It is classical to represent finitely generated subgroups of $F$ by
pointed $A$-graphs.  Let $H$ be a finitely generated subgroup of $F$
(we write $H\lefg F$) and let $h_1,\ldots,h_m \in F(A)$ be a set of
generators of $H$.  Let $\Gamma_0(h_1,\ldots,h_m)$ be the dual
$A$-graph which consists of a bouquet of $m$ loops, labeled by the
$h_i$, around a distinguished vertex, usually denoted by 1.  We denote
by $\Gamma(H)$ the reduced $A$-graph $\red(\Gamma_0(h_1,\ldots,h_m))$.
Observe that this construction coincides with the application of the
so-called Stallings foldings \cite{Stallings,KM}: it is well-known
that the pair $(\Gamma(H),1)$ depends on $H$ only, not on the choice
of the generating set $\{h_1,\ldots,h_m\}$, and we call it
\textit{the} (\textit{graphical}) \textit{representation of $H$}.
$\Gamma(H)$ is a connected reduced $A$-graph, in which no vertex
different from 1 is an endpoint.

Conversely, if $\Gamma$ is a connected reduced $A$-graph, 1 is a
vertex of $\Gamma$ and $\Gamma$ has no endpoint except maybe for 1,
there exists a unique subgroup $H \lefg F$ such that $(\Gamma,1)$ is
the representation of $H$.  In that situation, let $T$ be a spanning
tree of the $A$-graph $\Gamma$, and for each vertex $x$, let $u_{x}$
be the only reduced word labeling a path from $1$ to $x$ inside the
tree $T$.  For each positively labeled edge $e = (x,a,y)$ (that is,
$a\in A$), let $h_{e} = u_{x}a\bar u_{y}$.  Then a basis of $H$
consists of the elements $h_{e}$, where $e$ runs over the positively
labeled edges not in $T$ \cite{Stallings,KM}.

\begin{example}\label{example subgroups}
    Let $H_1 = \langle a^2b\inv, b^2a\inv\rangle$.  With reference to
    the graphs in Figure~\ref{fig graphs}, we see that $\Gamma_0 =
    \Gamma_0(a^2b\inv, b^2a\inv)$ (with distinguished vertex 1) and
    $(\Gamma_1,1)$ is the graphical representation of $H_1$.  Let $H_2
    = \langle a^3b\inv a\inv, ab^2a^{-2}\rangle$ and $H_3 = \langle
    b\inv ab\inv ab, b\inv a\inv b^3\rangle$.  The graphical
    representation of $H_2$ is $(\Gamma_2, 2)$ and the graphical
    representation of $H_3$ is $(\Gamma_3,3)$.
\end{example}

\begin{fact}\label{fact flower and fold}
    In view of Fact~\ref{fact reducing and trimming}, and if
    $\sum_{i=1}^m|h_i| = n$, computing $\Gamma(H)$ takes time
    $O(n^2\log(nr))$.  Given an $n$-vertex reduced $A$-graph $\Gamma$
    and a vertex 1, and letting $H$ be the subgroup represented by
    $(\Gamma,1)$, one can compute in time $O(n^2\log(nr))$ a basis of
    $H$, whose elements are words of length at most $2n$ \cite[Sec.
    3.3]{SilvaWeil}.
    
    If $H = \langle u\rangle$ is a cyclic subgroup of $F$, generated
    by a word of length $n\ne 0$, then $\{u\}$ is a basis of $H$ and
    $\Gamma(\langle u\rangle) = \red(\Gamma(u))$ is computed in time
    $O(n\log r)$ by Fact~\ref{fact reducing and trimming}.
\end{fact}    

Let $H\lefg F$, with representation $(\Gamma(H),1)$.  Let us say that
\textit{$H$ is cyclically reduced} if 1 is not an endpoint, that is,
if $\Gamma(H)$ is cyclically reduced.  In the general case, let $K$ be
the subgroup represented by $(\cc(\Gamma(H)),\beta(1))$.  It is
well-known that $K = \chi_{b(1)}(H)$, where $\chi_u$ is the
conjugation $x \mapsto u\inv xu$.  It follows that, for each subgroup
$H'\lefg F$, $H'$ is a conjugate of $H$ if and only if
$\cc(\Gamma(H')) = \cc(\Gamma(H))$.  As a consequence, we say that
$\cc(\Gamma(H))$ is \textit{the} (\textit{graphical})
\textit{representation of the conjugacy class} $[H]$.

We note that if $H = \langle u\rangle$ is a cyclic subgroup, then the
subgroup $H$ is cyclically reduced if and only if the word $u$ is
cyclically reduced.

In the sequel, the size of $\Gamma(H)$ and of $\cc(\Gamma(H))$ will be
taken to be measures of the \textit{size} of $H$ and $[H]$, and we
will write $|H| = |\Gamma(H)|$ and $|[H]| = |\cc(\Gamma(H))|$.  In
particular, $|[H]|$ is the minimum size of subgroups in the
conjugacy class $[H]$, and it is equal to the size of any
cyclically reduced conjugate of $H$.

\begin{example}\label{example cyclic core}
    Let $H_1$, $H_2$ and $H_3$ be the subgroups discussed in
    Example~\ref{example subgroups}.  In view of Figure~\ref{fig
    graphs}, we see that $H_1$ is cyclically reduced, and that $H_2 =
    \chi_{a\inv}(H_1)$ and $H_3 = \chi_{ab}(H_1)$.
\end{example}

%%%%%%%%%%%%%%%%%%%%%%%
\subsection{Action of an automorphism}\label{sec action}

Let $\phi\in\Aut(F)$ and let $\Gamma = (V,E)$ be a reduced $A$-graph 
with a designated vertex 1.  Let
$\phi_\bullet(\Gamma)$ be the $A$-graph obtained after the following steps:
\begin{itemize}
    \item[(a)] replace each $a$-labeled edge by a path labeled by the
    word $\phi(a)$ with the same endpoints: if $\phi(a) = a_1\cdots
    a_m$ ($a_i\in\tilde A$) and $(x_0,a,x_m)\in E$, remove the edges
    $(x_0,a,x_m)$ and $(x_m,\bar a,x_0)$, add new vertices
    $x_1,\ldots,x_{m-1}$ and add edges $(x_{i-1},a_i,x_i)$ and
    $(x_i,\bar a_i,x_{i-1})$ for each $1\le i \le m$;
    
    \item[(b)] reduce the resulting $A$-graph;
    
    \item[(c)] repeatedly trim all the endpoints different from 1.
\end{itemize}
If $\Gamma$ is cyclically reduced, we also denote by $\phi(\Gamma)$
the cyclically reduced graph $\cc(\phi_\bullet(\Gamma))$.

\begin{figure}[ht]
    \centering
    \begin{picture}(126,31)(0,-31)
	
	\node[Nw=5.0,Nh=5.0,Nmr=1.0](n22)(0.0,-1.0){1}

	\node[Nw=5.0,Nh=5.0,Nmr=1.0](n23)(22.0,-1.0){}

	\node[Nw=5.0,Nh=5.0,Nmr=1.0](n24)(11.0,-17.0){}

	\put(11.0,-28.0){$\Gamma_1$}

	\node[Nw=5.0,Nh=5.0,Nmr=1.0](n2)(39.0,-1.0){1}

	\node[Nw=5.0,Nh=5.0,Nmr=1.0](n3)(61.0,-1.0){}

	\node[Nw=5.0,Nh=5.0,Nmr=1.0](n4)(81.0,-1.0){}

	\node[Nw=5.0,Nh=5.0,Nmr=1.0](n5)(61.0,-17.0){}

	\node[Nw=5.0,Nh=5.0,Nmr=1.0](n6)(81.0,-17.0){}
	
	\put(61.0,-28.0){$\phi_\bullet(\Gamma_1)$}

	\node[Nw=5.0,Nh=5.0,Nmr=1.0](n30)(102.0,-1.0){}

	\node[Nw=5.0,Nh=5.0,Nmr=1.0](n40)(124.0,-1.0){}

	\node[Nw=5.0,Nh=5.0,Nmr=1.0](n50)(102.0,-17.0){}

	\node[Nw=5.0,Nh=5.0,Nmr=1.0](n60)(124.0,-17.0){}
	
	\put(102.0,-28.0){$\phi(\Gamma_1)$}

    \drawedge(n2,n3){$b$}

    \drawedge[ELside=r](n3,n5){$a$}

    \drawedge[ELside=r](n4,n3){$a$}

    \drawedge(n5,n6){$a$}

    \drawedge(n4,n6){$b$}

    \drawloop[loopdiam=6.0,loopangle=180.0](n5){$b$}

    \drawedge[ELside=r](n30,n50){$a$}

    \drawedge[ELside=r](n40,n30){$a$}

    \drawedge(n50,n60){$a$}

    \drawedge(n40,n60){$b$}

    \drawloop[loopdiam=6.0,loopangle=180.0](n50){$b$}

    \drawedge(n22,n23){$a$}

    \drawedge[ELside=r](n22,n24){$b$}

    \drawedge[curvedepth=3.0](n23,n24){$a$}

    \drawedge[curvedepth=3.0](n24,n23){$b$}

\end{picture}
\caption{Action of an automorphism on an $A$-graph}
\label{fig automorphic image}
\end{figure}
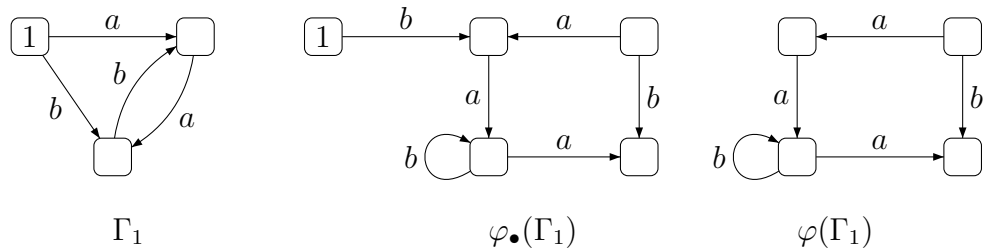

\begin{example}\label{example automorphic image}
    Let $\Gamma_1$ be as in Example~\ref{example graphs} and let
    $\phi\in\Aut(F)$ be given by $\phi(a) = ba\inv$ and $\phi(b) =
    bab\inv$.  The graphs $\Gamma_1$, $\phi_\bullet(\Gamma_1)$ and
    $\phi(\Gamma_1)$ are shown in Figure~\ref{fig automorphic image}.
\end{example}

\begin{fact}\label{fact image by phi}
    Let $\ell$ be the maximum length of the image of a letter by an
    automorphism $\phi$ (so $\phi$ is given by a $r$-tuple of words of
    length at most $\ell$) and let $\Gamma$ be an $n$-vertex reduced
    $A$-graph.  In view of Fact~\ref{fact reducing and trimming}, the
    complexity of computing $\phi_\bullet(\Gamma)$ is
    $O(n^2\ell^2r^2\log(n\ell r))$.
    
    If $\Gamma$ is a cyclic word, $\Gamma = \Gamma(u)$ with $u\in
    F(A)$ cyclically reduced and $|u| = n$, then $\phi_\bullet(\Gamma)
    = \red(\Gamma(\phi(u)))$, which is computed in time $O(n\ell\log
    r)$.
\end{fact}

It is easy to verify that if $\phi\in\Aut(F)$ and $H\lefg F$, then
$\phi(H)$ is represented by $(\phi_\bullet(\Gamma(H)),1)$.  If in
addition $H$ is cyclically reduced, the conjugacy class
$\phi([H])$ is represented by $\phi(\Gamma(H))$.

The \textit{Whitehead minimization problem} (WMP) \textit{for finitely
generated subgroups} (resp.  \textit{for conjugacy classes of finitely
generated subgroups}) consists in finding a minimum size element $X'$
of the automorphic orbit of a given finitely generated subgroup (resp.
conjugacy class of finitely generated subgroups) $X$, and an
automorphism $\phi$ such that $X' = \phi(X)$.

If the input of the WMP is a conjugacy class of cyclic subgroups, that
is, a cyclic word, we talk of the \textit{Whitehead minimization
problem for cyclic words}.  The \textit{Whitehead minimization problem
for words} consists in finding a minimum length element $u'$ of the
automorphic orbit of a given reduced word $u$, and an automorphism
$\phi$ such that $u' = \phi(u)$.
    
We will see in Section~\ref{sec W theorem} that all these problems
reduce to the problem for conjugacy classes of finitely generated
subgroups.

%%%%%%%%%%%%%%%%%%%%%%
\subsection{Gersten's theorem and the Whitehead method}\label{sec W theorem}

It is well-known that the group $\Aut(F)$ of automorphisms of $F$ is
finitely generated.  One particular finite set of generators of
$\Aut(F)$ is the set of so-called Whitehead automorphisms (relative to
the choice of the basis $A$), whose precise definition will be given
in Section~\ref{sec W automorphisms} below.

The first key element of the algorithm presented here, is the
following statement, due to Gersten \cite[Corol. 2]{Gersten}.

\begin{theorem}\label{W theorem}
    Let $H$ be a cyclically reduced subgroup of $F(A)$.  If there
    exists an automorphism $\phi$ of $F$ such that $|\phi([H])| <
    |[H]|$, then there exists such an automorphism among the
    Whitehead automorphisms.
\end{theorem}

\begin{remark}
    This theorem is a generalization of a fundamental result of
    Whitehead (\cite{Whitehead}, see \cite[Sec.  I.4]{LS}), which
    concerns the cyclic word case --- the case where $H = \langle
    u\rangle$ for some cyclically reduced word $u$.
\end{remark}    

This implies the following algorithm --- the so-called
\textit{Whitehead method} --- to solve the WMP for a conjugacy class
of finitely generated subgroups $[H]$.  The input of the algorithm
is a cyclically reduced subgroup $H$, or rather the cyclically reduced
$A$-graph $\Gamma(H)$.  The output of the algorithm is a cyclically
reduced subgroup $H'$ and a tuple $\vec\phi = (\phi_m,\ldots,\phi_1)$
of Whitehead automorphisms, such that $[H'] = \phi_m \circ \cdots
\circ \phi_1([H])$ and $[H']$ has minimum size in the
automorphic orbit of $[H]$.

First let $\vec\phi = (\id)$ and $\Gamma = \Gamma(H)$.  Then
repeatedly apply the following steps: try every Whitehead automorphism
$\psi$ until $|\psi(\Gamma)| < |\Gamma|$; if such a $\psi$ exists,
replace $\Gamma$ by $\psi(\Gamma)$ and $\vec\phi$ by $(\psi,
\vec\phi)$; otherwise, stop and output $\Gamma$ and $\vec\phi$.  At
each step, the size of $\Gamma$ decreases by at least one unit, so
this procedure terminates after at most $|\Gamma|$ iterations.
Finally, in order to output a basis of (a possible value of) $H'$,
choose arbitrarily a vertex 1 in $\Gamma$ and use the procedure
discussed in Section~\ref{sec subgroups}.

Let us give an estimate of the complexity of this algorithm.  

\begin{fact}\label{fact WMP conjclass}
    We first note that the cost of the construction of $\Gamma(H)$, if
    the input is a set of generators of $H$ of total length $n$, is
    $O(n^2\log(nr))$ (Fact~\ref{fact flower and fold}).  Moreover,
    $\Gamma(H)$ has at most $n$ vertices.
    
    As we will see, a Whitehead automorphism maps every letter to a
    word of length at most 3, so finding the image of a cyclically
    reduced graph of size $n$ under a Whitehead automorphism also
    takes time $O(n^2r^2\log(nr))$ (Fact~\ref{fact image by phi}).
    
    Thus, if $f(r)$ is the cardinality of the set of Whitehead
    automorphisms (of a rank $r$ free group), each iterating step of
    the algorithm may require $f(r)$ steps, each of which consists in
    computing the image of a cyclically reduced graph of length at
    most $n$ under a Whitehead automorphism, and hence has complexity
    $O(n^2r^2\log(nr))$.  There are at most $n$ iterating steps, so
    the iterating part of the algorithm has time complexity
    $O(n^3r^2f(r)\log(nr))$.
    
    Finally, retrieving a basis of $H'$ from the ultimate value of
    $\Gamma$ takes time $O(n^2\log(nr))$ (Fact~\ref{fact flower and
    fold}).  That is negligible in front of $n^3r^2f(r)\log(nr)$, so
    the total complexity of the algorithm is $O(n^3r^2f(r)\log(nr))$,
    which is polynomial in $n$ and exponential in $r$ as we shall see
    in Section~\ref{sec W automorphisms}.
    
    In the cyclic word case, that is, the case where $H = \langle
    u\rangle$ with $u\in F(A)$ cyclically reduced and $|u| = n$, the
    complexity is $O(n^2f(r)\log r)$.
\end{fact}

\begin{remark}
    Let $(H',\vec\phi)$ be the output of the algorithm on input $H$.
    The complexity discussion above shows that $\vec\phi =
    (\phi_m,\ldots,\phi_1)$ with $m < n$.  As we will see in
    Section~\ref{sec W automorphisms}, the length of the image of a
    letter in a Whitehead automorphism is at most 3, so the length of
    $\phi_m \circ \cdots \circ \phi_1(a)$ ($a\in A$) may be
    exponential in $m$, and the computation of $\phi_m \circ \cdots
    \circ \phi_1$ may take time exponential in $m$.  This possible
    exponential explosion is the reason why we choose to output a
    tuple of Whitehead automorphisms rather than their composition.
    
    An easy example for this exponential explosion is provided by
    the (Whitehead) automorphisms $\alpha\colon a\mapsto ab;\ b\mapsto
    b$ and $\beta\colon a\mapsto a;\ b\mapsto ba$.  Then
    $\beta\circ\alpha\colon a\mapsto aba;\ b\mapsto ba$, and the
    length of $(\beta\circ\alpha)^n(a)$ (resp.
    $(\beta\circ\alpha)^n(b)$) is the sum of the entries of the first
    (resp.  second) column of ${\scriptstyle\pmatrix{\scriptstyle 2 & \scriptstyle
    1 \cr \scriptstyle 1 & \scriptstyle 1}}^n$.  It is well-known that
    the asymptotic behavior of these numbers as $n$ tends to infinity,
    is $O(\rho^n)$ where $\rho$ is the dominant eigenvalue of that
    matrix, namely, $\rho = (3 + \sqrt 5)/2$.
\end{remark}    

The above algorithm can be modified to solve the WMP for a finitely
generated subgroup $H$ in $F(A)$.  A minimum size element $K$ of the
automorphic orbit of $H$ is in particular a minimum size element of
its own conjugacy class, and hence $K$ is cyclically reduced.
Moreover, the class $[K]$ must be a minimum size element in the
automorphic orbit of $[H]$.  Thus a minimum size element of the orbit
of $H$ is obtained by computing a minimum size element of the orbit of
$[H]$, say $[K]$, and choosing arbitrarily a cyclically reduced
subgroup in $[K]$.

More precisely, the algorithm is as follows.  First let
$(\Gamma(H),1)$ be the graphical representation of $H$, let $\Gamma =
\cc(\Gamma(H))$ and $\vec\phi = (\chi_{b(1)})$.  We note that $\Gamma
= {\chi_{b(1)}}_\bullet(\Gamma(H))$.  Next, rename as 1 the vertex
$\beta(1)$ of $\Gamma$.  Then repeatedly apply the following steps:
try every Whitehead automorphism $\psi$ until $|\psi(\Gamma)| <
|\Gamma|$; if such a $\psi$ exists, consider the pointed $A$-graph
$(\psi_\bullet(\Gamma),1)$, replace $\vec\phi$ by
$(\chi_{b(1)}\circ\psi, \vec\phi)$, rename $\beta(1)$ as 1, and
replace $\Gamma$ by $\psi(\Gamma) = \cc(\phi_\bullet(\Gamma))$;
otherwise, stop and output $\Gamma$ and $\vec\phi$.  Finally, construct a
basis of the subgroup $H'$ represented by $(\Gamma,1)$.

Finally, this last algorithm can be used to solve the WMP for words: a
minimum length element in the automorphic orbit of a word $u\in F(A)$
is necessarily a cyclically reduced word $u'$ such that $\langle
u'\rangle$ is a solution of the WMP on input $\langle u\rangle$.
Therefore, it suffices to apply the above algorithm on input $\langle
u\rangle$, letting $u' = u$ at the beginning of the algorithm, and
updating $u'$ to $\chi_{b(1)}\circ\psi(u')$ at each iterating step.
We note that the length of $u'$ never exceeds $|\Gamma|$, and hence
never exceeds $|u| = n$.

\begin{fact}\label{fact WMP subgroup}
    The extra work required by this modified algorithm (see
    Fact~\ref{fact WMP conjclass}), namely to compute $\cc(\Gamma(H))$
    and to compose at most $n$ Whitehead automorphisms with
    conjugations by words of length at most $n$, takes time
    $O(n^2r\log r)$, which is negligible in front of
    $n^3r^2f(r)\log(nr)$.  So the time complexity of this algorithm is
    again $O(n^3r^2f(r)\log(nr))$.
    
    In the cyclic subgroup case, as well as in the word case, the
    complexity is $O(n^2f(r)\log r)$.
\end{fact}

%%%%%%%%%%%%%%%%%%%%%%
\subsection{Whitehead automorphisms}\label{sec W automorphisms}

We now review the definition of the \textit{Whitehead automorphisms}
of $F$, relative to the choice of the basis $A$, see for instance
\cite[Sec.  I.4]{LS}.

There are two kinds of Whitehead automorphisms.  The \textit{first
kind} consists of the automorphisms that permute the set $\tilde A$.
We observe that these are exactly the length-preserving automorphisms
of $F(A)$, that is, the automorphisms $\phi$ such that $|\phi(u)| =
|u|$ for each $u\in F(A)$.  Each is uniquely specified by a
permutation $\sigma$ of $A$ and an $A$-tuple $\mathbf{x} = (x_a)_{a\in
A} \in \{\pm1\}^A$: the automorphism specified by $\sigma$ and
$\mathbf{x}$ maps each letter $a$ to $\sigma(a)^{x_a}$.  In
particular, there are $r!\,2^r$ length-preserving (Whitehead)
automorphisms.

Let $v\in \tilde A$.  We define a \textit{$v$-cut of $\tilde A$} to be
a subset $Y\subseteq \tilde A$ containing $v$ and avoiding $\bar v$.
Each pair $(v,Y)$ of a letter $v\in \tilde A$ and a $v$-cut $Y$
defines a \textit{Whitehead automorphism of the second kind} $\phi$ as
follows: $\phi(v) = v$ and for each $a\in A\setminus\{v,\bar v\}$,
    $$\phi(a) = v^\gamma a v^\rho \textrm{ where }
    \gamma = \cases{-1 & if $\bar a\in Y$, \cr 0 & otherwise;}
    \qquad \rho = \cases{1 & if $a\in Y$, \cr 0 & otherwise.}$$

\begin{remark}\label{remark gamma rho}
    By inverting both sides of this formula, which specifies the
    images of the letters in $A$ under $\phi$, we find that the same
    formula also holds for the letters in $\bar A$: if $a\in \bar A$
    (and $\bar a\in A$) and $a\ne v,\bar v$, then $\phi(a) = \phi(\bar
    a)\inv = v^\gamma a v^\rho$ where $\gamma = -1$ if $\bar a\in Y$
    and $\rho = 1$ if $a\in Y$.
\end{remark}

Observe that, if $Y$ is reduced to the singleton $\{v\}$, then the
resulting Whitehead automorphism is the identity.  Apart from this
particular case, no Whitehead automorphism of the second kind is
length-preserving, and the automorphisms specified by different pairs
$(v,Y)$ and $(v',Y')$ are distinct.  In particular, if $\W(A)$ denotes
the set of non-identity Whitehead automorphisms of the second kind,
then $\card(\W(A)) = 2r\,(2^{2r-2}-1) = r\,(2^{2r-1}-2)$.

Finally, we note that in the algorithms to solve the WMP discussed in
Section~\ref{sec W theorem}, the set of all Whitehead automorphisms
can be replaced throughout by $\W(A)$, since we care only about the
length of words and cyclic words, and since $\W(A)$ is preserved by
composition with the length-preserving Whitehead automorphisms.  That
is, the function $f(r)$ is Facts~\ref{fact WMP conjclass} and
\ref{fact WMP subgroup} can be taken equal to $r2^{2r}$.  In
particular, we have the following fact.

\begin{fact}
    The algorithms given in Section~\ref{sec W theorem} to solve the
    WMP for conjugacy classes of subgroups or for subgroups, take time
    $O(n^3 r^3 4^r\log(nr))$.
    
    The algorithms to solve the WMP for words, cyclic words and cyclic
    subgroups take time $O(n^2 r\, 4^r\log r)$.
\end{fact}

%%%%%%%%%%%%%%%%%%%%%%
\section{Choice of a best Whitehead automorphism}\label{sec choice}

The algorithms given above to solve the WMP are exponential in the
rank $r$ of $F$ because every element of $\W(A)$ may have to be tested
at each iteration of the algorithm.  Our point is to show that, given
a cyclically reduced $A$-graph $\Gamma$, one can in polynomial time
(in $r$ and in $|\Gamma|$) find an element $\phi$ of $\W(A)$ that
minimizes $|\phi(\Gamma)|$ --- thereby determining in particular
whether $|\Gamma|$ is minimal.  Our first tool for this purpose is a
generalization of the classical Whitehead graph associated with a
cyclic word.

%%%%%%%%%%%%%%%%%%%%%%
\subsection{Whitehead hypergraph of a cyclically reduced $A$-graph}\label{sec hypergraph}

A \textit{hypergraph} (here, an undirected one) is a triple $G =
(B,D,\kappa)$ where $B$ and $D$ are sets and $\kappa\colon D
\rightarrow \calP(B)$ is such that $\kappa(d)\ne\emptyset$ for each
$d\in D$.  The elements of $B$ are called \textit{vertices}, the
elements of $D$ are called \textit{hyperedges}.  It is understood that
there can be several hyperedges with the same $\kappa$-image.  A
hypergraph in which every hyperedge has cardinality at most (resp.
exactly) 2 can be identified naturally with an undirected graph (resp.
a loop-free undirected graph).

Let $\Gamma = (V,E)$ be a cyclically reduced $A$-graph.  The
\textit{Whitehead hypergraph} of $\Gamma$, written $W_\Gamma$, is
defined as follows.  Its vertex set is $\tilde A$.  Its hyperedge set
$D$ is in bijection with $V$, $D = \{d_x \mid x\in V\}$, and for each
vertex $x$ of $\Gamma$, $\kappa(d_x)$ is the hyperlink of $x$,
$\kappa(d_x) = \hl(x) = \{\lambda(e) \mid e\in E,\ \omega(e) = x\}$.

Note that every hyperedge of $W_\Gamma$ has a $\kappa$-image with at 
least two elements since $\Gamma$ is cyclically reduced. 

\begin{example}\label{ex Whitehead graph}
    Let $\Gamma$ be the $A$-graph shown below, where $A =
    \{a,b,c,d,e\}$.  Since the vertex set of $\Gamma$ is
    $\{1,2,3,4,5,6\}$, the Whitehead hypergraph $W_\Gamma$ has vertex
    set $\tilde A$ and hyperedge set $\{d_i \mid 1\le i \le 6\}$, with

    \begin{picture}(65,45)(-20,-40)

	\node[Nw=5.0,Nh=5.0,Nmr=2.0](n1)(2.0,-0.0){1}

	\node[Nw=5.0,Nh=5.0,Nmr=2.0](n2)(2.0, -20.0){2}

	\node[Nw=5.0,Nh=5.0,Nmr=2.0](n3)(2.0,-40.0){3}

	\node[Nw=5.0,Nh=5.0,Nmr=2.0](n4)(25.0,-0.0){4}

	\node[Nw=5.0,Nh=5.0,Nmr=2.0](n5)(25.0,-20.0){5}

	\node[Nw=5.0,Nh=5.0,Nmr=2.0](n6)(25.0,-40.0){6}

	\drawedge(n1,n4){$a$}

	\drawedge[ELside=r](n1,n2){$b$}

	\drawedge[ELside=r](n2,n3){$d$}

	\drawedge(n5,n3){$c$}

	\drawedge(n5,n6){$d$}

	\drawedge(n3,n6){$e$}

	\drawedge(n4,n5){$a$}

	\drawedge[curvedepth=3.0](n2,n5){$c$}

	\drawedge[curvedepth=3.0](n5,n2){$a$}
	
	\put(50,-2){$\kappa(d_1) = \{\bar a, \bar b\}$}

	\put(50,-9){$\kappa(d_2) = \{a,b,\bar c,\bar d\}$}

	\put(50,-16){$\kappa(d_3) = \{c,d,\bar e\}$}

	\put(50,-23){$\kappa(d_4) = \{a,\bar a\}$}

	\put(50,-30){$\kappa(d_5) = \{a,\bar a,c,\bar c, \bar d\}$}

	\put(50,-37){$\kappa(d_6) = \{d,e\}$}

    \end{picture}
\end{example}
    
\begin{example}
    Let $u$ be a cyclically reduced word.  Then the hypergraph
    $W_{\Gamma(\langle u\rangle)}$ is in fact a graph ($\kappa$ maps
    each hyperedge to a pair of distinct vertices), denoted by $W_u$,
    which coincides with the classical notion of the \textit{Whitehead
    graph of a cyclically reduced word} \cite[Sec.  I.7]{LS}.
\end{example}

\begin{fact}\label{fact complexity Wh graph}
    If $|\Gamma| = n$, then $W_\Gamma$ has $2r$ vertices, $n$
    hyperedges and can be constructed in time $O(nr\log r)$.
\end{fact}

%%%%%%%%%%%%%%%%%%%%%%
\subsection{Applying a Whitehead automorphism}\label{sec technical}

We now come to the technical core of this paper: given a cyclically
reduced $A$-graph $\Gamma$ and a Whitehead automorphism $\phi\in
\W(A)$, specified by a pair $(v,Y)$, we give an exact formula for the
size difference $|\phi(\Gamma)| - |\Gamma|$.  In fact, this formula is
already known: it was established by Gersten \cite[Prop.  1]{Gersten},
proved in Kalajd\v zievski \cite[Prop.  10.3]{Kalajd}, and it is a
generalization of a result of Whitehead \cite[Prop.  I.4.16]{LS},
which covers the cyclic word case.  Our contribution here consists in
rewording it in graph-theoretic terms, and possibly in a clearer
demonstration.  We note that this formula is an essential ingredient
in the proof of Whitehead's theorem and its generalization by Gersten
(Theorem~\ref{W theorem} above), see \cite{LS,Gersten,Kalajd}.

Let $G$ be an undirected hypergraph, with vertex set $V(G)$.  We
define the \textit{capacity} of a subset $Y \subseteq V(G)$ to be the
number $\capac_G(Y)$ of hyperedges $e$ of $G$ such that $\kappa(e)$
meets both $Y$ and its complement $Y^c$.  If $v\in V(G)$, the
\textit{degree of $v$} is the number $\deg_G(v)$ of hyperedges whose
$\kappa$-image contains $v$ (that is, that are adjacent to $v$).  We
show the following result.

\begin{proposition}\label{prop technical}
    Let $\Gamma$ be a cyclically reduced $A$-graph, let $v\in\tilde A$
    and let $Y\subseteq \tilde A$ be a $v$-cut of $\tilde A$ (i.e., a
    set containing $v$ and not $\bar v$).  Let $\phi$ be the Whitehead
    automorphism specified by the pair $(v,Y)$.  Then we have
    $$|\phi(\Gamma)| - |\Gamma| = \capac_{W_\Gamma}(Y) -
    \deg_{W_\Gamma}(v).$$
\end{proposition}

\proof
Let $\Gamma =(V,E)$, $v$, $Y$ and $\phi$ be as in the statement.  We
first examine in detail the construction of $\phi(\Gamma)$.  The first
step is to construct the $A$-graph $\Gamma' = (V',E')$ as follows.
For each vertex $x\in V$, if $\hl_{\Gamma}(x) \cap (Y\setminus\{v\})
\neq \emptyset$, we let $u_x$ be a new vertex.  If $\hl_{\Gamma}(x)
\subseteq Y^c\cup\{v\}$, $u_x$ is undefined.  We let $V' = V \cup
\{u_x \mid \textrm{ $u_x$ exists}\}$.

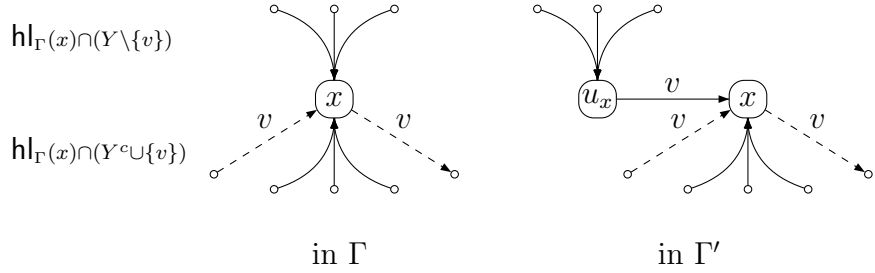
\begin{figure}[ht]
    \centering
\begin{picture}(114,36)(0,-36)
    
    \put(0,-7){$\scriptstyle \hl_\Gamma(x) \cap (Y \setminus\{v\})$}

    \node[Nw=1.0,Nh=1.0](n1)(27.0,-24.0){}

    \node[Nw=1.0,Nh=1.0](n2)(35.0,-2.0){}

    \node[Nw=1.0,Nh=1.0](n3)(43.0,-2.0){}

    \node[Nw=1.0,Nh=1.0](n4)(51.0,-2.0){}

    \put(0,-22){$\scriptstyle \hl_\Gamma(x) \cap (Y^c \cup\{v\})$}

    \node[Nw=1.0,Nh=1.0](n5)(35.0,-26.0){}

    \node[Nw=1.0,Nh=1.0](n6)(43.0,-26.0){}

    \node[Nw=1.0,Nh=1.0](n7)(51.0,-26.0){}

    \node[Nw=1.0,Nh=1.0](n8)(59.0,-24.0){}

    \node[Nw=5.0,Nh=5.0,Nmr=2.0](n9)(43.0,-14.0){$x$}

    \drawedge[dash={1.0 1.0}{0.0}](n1,n9){$v$}

    \drawedge[dash={1.0 1.0}{0.0}](n9,n8){$v$}

    \drawedge[curvedepth=3.0](n2,n9){}

    \drawedge(n3,n9){}

    \drawedge[curvedepth=-3.0](n4,n9){}

    \drawedge[curvedepth=-3.0](n5,n9){}

    \drawedge(n6,n9){}

    \drawedge[curvedepth=3.0](n7,n9){}

    \put(40,-36){in $\Gamma$}

    \node[Nw=1.0,Nh=1.0](n11)(82.0,-24.0){}

    \node[Nw=1.0,Nh=1.0](n12)(70.0,-2.0){}

    \node[Nw=1.0,Nh=1.0](n13)(78.0,-2.0){}

    \node[Nw=1.0,Nh=1.0](n14)(86.0,-2.0){}

    \node[Nw=1.0,Nh=1.0](n15)(90.0,-26.0){}

    \node[Nw=1.0,Nh=1.0](n16)(98.0,-26.0){}

    \node[Nw=1.0,Nh=1.0](n17)(106.0,-26.0){}

    \node[Nw=1.0,Nh=1.0](n18)(114.0,-24.0){}

    \node[Nw=5.0,Nh=5.0,Nmr=2.0](n19)(98.0,-14.0){$x$}

    \node[Nw=5.0,Nh=5.0,Nmr=2.0](n20)(78.0,-14.0){$u_x$}

    \drawedge[dash={1.0 1.0}{0.0}](n11,n19){$v$}

    \drawedge[dash={1.0 1.0}{0.0}](n19,n18){$v$}

    \drawedge[curvedepth=3.0](n12,n20){}

    \drawedge(n13,n20){}

    \drawedge[curvedepth=-3.0](n14,n20){}

    \drawedge[curvedepth=-3.0](n15,n19){}

    \drawedge(n16,n19){}

    \drawedge[curvedepth=3.0](n17,n19){}

    \drawedge(n20,n19){$v$}

    \put(86,-36){in $\Gamma'$}

\end{picture}
\caption{From $\Gamma$ to $\Gamma'$}
\label{fig gamma'}
\end{figure}

The set $E'$ consists of the following edges.  All the $v$- and $\bar
v$-labeled edges of $\Gamma$ are also in $E'$.  Next, if $x\in V$ and
$u_x$ exists, then there is a $v$-labeled edge from $u_x$ to $x$ and a
$\bar v$-labeled edge from $x$ to $u_x$.  Finally, for each
$a$-labeled ($a\ne v,\bar v$) edge from $x$ to $y$ in $\Gamma$, there
is an $a$-labeled edge in $\Gamma'$
\begin{itemize}
    \item from $u_x$ to $u_y$ if $a,\bar a\in Y$,
    
    \item from $x$ to $u_y$ if $a\in Y$, $\bar a\not\in Y$,
    
    \item from $u_x$ to $y$ if $a\not\in Y$, $\bar a\in Y$,
    
    \item from $x$ to $y$ if $a,\bar a\not\in Y$.
\end{itemize}
The transformation of $\Gamma$ into $\Gamma'$ is local in the
following sense.  Around each vertex $x$, we separate
$\hl_{\Gamma}(x)$ into $\hl_{\Gamma}(x)\cap (Y\setminus \{v\})$ and
$\hl_{\Gamma}(x)\cap (Y^c \cup\{ v\})$ and, when the first set is
non-empty, we push this fragment of $\hl_\Gamma(x)$ away from $x$ by
introducing a new $v$-labeled edge (see Figure~\ref{fig gamma'}).

An observation that will be important in the sequel is that, by
construction, the vertices $x\in V$ satisfy $\hl_{\Gamma'}(x)\subseteq
Y^c \cup \{ v\}$, while the new vertices of the form $u_x$ satisfy
$\hl_{\Gamma'}(u_x) \subseteq (Y\setminus \{ v\})\cup \{ \bar v\}$.

We note that for each $a$-labeled edge from $x$ to $y$ in $\Gamma$, we
now have a path in $\Gamma'$ from $x$ to $y$, labeled by the word
$\phi(a)$, and that $\Gamma'$ consists of the collection of these
paths.  Thus $\phi(\Gamma)$ is obtained by first reducing $\Gamma'$,
and then taking the cyclic core of the resulting $A$-graph, see
Section~\ref{sec action}.

We now consider whether $\Gamma'$ is reduced.  Let $x\in V$ be such
that $u_x$ exists.  Then $\link_{\Gamma'}(u_x)$ consists of a $\bar
v$-labeled edge, and a non-empty set in bijection with the set of
edges in $\link_\Gamma(x)$ with a label in $Y\setminus \{v\}$.  In
particular, the labeling map $\lambda'$ is injective on
$\link_{\Gamma'}(u_x)$.  Moreover, $u_x$ is not an endpoint.

Now let $x$ be a vertex of $\Gamma'$, in $V$. There are 3 cases. We let
\begin{itemize}
    \item $V_1$ be the set of $x\in V$ such that $u_x$ does not exist,
    that is, $\hl_\Gamma(x) \subseteq Y^c \cup\{v\}$;
    
    \item $V_2$ be the set of $x\in V$ such that $u_x$ exists and
    $\link_\Gamma(x)$ contains no $v$-labeled edge, that is,
    $\hl_{\Gamma}(x) \cap (Y\setminus\{v\}) \neq \emptyset$ and
    $v\not\in\hl_\Gamma(x)$;

    \item $V_3$ be the set of $x\in V$ such that $u_x$ exists and
    $\link_\Gamma(x)$ contains a $v$-labeled edge, that is,
    $\hl_{\Gamma}(x) \cap (Y\setminus\{v\}) \neq \emptyset$ and
    $v\in\hl_\Gamma(x)$.
\end{itemize}

\textbf{Case 1: $x\in V_1$}.  \enspace Then $\link_{\Gamma'}(x)$ is in
bijection with $\link_\Gamma(x)$, the labeling map $\lambda'$ is
injective on $\link_{\Gamma'}(x)$, and $x$ is not an endpoint.

\medskip

\textbf{Case 2: $x\in V_2$}.  \enspace Then $\link_{\Gamma'}(x)$
consists of a $v$-labeled edge plus a set in bijection with the subset
of $\link_\Gamma(x)$ of all edges labeled by letters in $Y^c$.  In
particular, the labeling map $\lambda'$ is injective on
$\link_{\Gamma'}(x)$.  Moreover, $x$ is an endpoint in $\Gamma'$ if
and only if $\hl_\Gamma(x)\subseteq Y$.

\medskip

\textbf{Case 3: $x\in V_3$}.  \enspace Then $\link_{\Gamma'}(x)$
consists of two $v$-labeled edges (one of them starting at $u_x$) and
a set in bijection with the subset of $\link_\Gamma(x)$ consisting of
the edges labeled by letters in $Y^c$.

\medskip

Thus $\Gamma'$ is non-reduced if and only if $V_3\ne\emptyset$, and in
that case, the first step in reducing $\Gamma'$ consists in performing
the elementary reductions that arise from the pairs of $v$-labeled
edges into the vertices $x\in V_3$.  We claim that the resulting
graph, say $\Gamma''$, is already reduced.

In order to justify this claim, let us consider the effect of such an
elementary reduction.  Since $x\in V_3$, $\Gamma$ has a $v$-labeled
edge from some $y\in V$ to $x$, and at least one $a$-labeled edge from
some $z\in V$ to $x$, with $a\in Y\setminus\{v\}$, see Figure~\ref{fig
xV3}.  In $\Gamma'$, there are $v$-labeled edges from $y$ and from
$u_x$ to $x$, and an $a$-labeled edge from $z'$ to $u_x$ (with
$z'\in\{z, u_z\}$).  By a previous observation,
$\hl_{\Gamma'}(y)\subseteq Y^c \cup \{ v\}$ and
$\hl_{\Gamma'}(u_x)\setminus \{ \bar v\} \subseteq Y\setminus \{ v\}$.
Thus, after the elementary reduction identifying $u_x$ to $y$, the
labeling function $\lambda''$ is injective on $\link_{\Gamma''}(y)$.
It follows that $\Gamma''$ is reduced.  Moreover $|\Gamma''|
=|\Gamma| - \card(V_2)$.

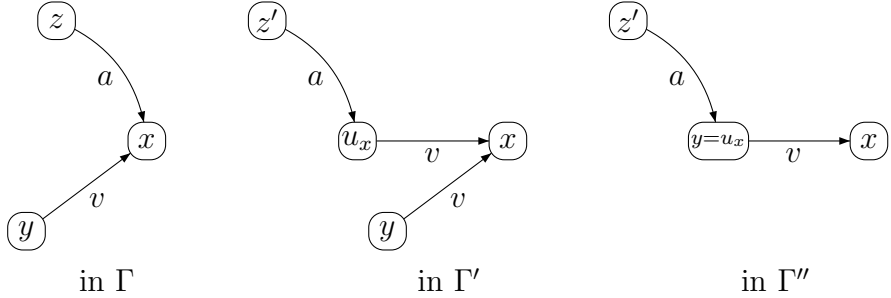
\begin{figure}[ht]
    \centering
\begin{picture}(112,38)(0,-38)

    \node[Nw=5.0,Nh=5.0,Nmr=2.0](n1)(4.0,-2.0){$z$}

    \node[Nw=5.0,Nh=5.0,Nmr=2.0](n2)(0.0,-30.0){$y$}

    \node[Nw=5.0,Nh=5.0,Nmr=2.0](n3)(16.0,-18.0){$x$}

    \drawedge[curvedepth=3,ELside=r](n1,n3){$a$}

    \drawedge[ELside=r](n2,n3){$v$}

    \put(7,-38){in $\Gamma$}

    \node[Nw=5.0,Nh=5.0,Nmr=2.0](n4)(32.0,-2.0){$z'$}

    \node[Nw=5.0,Nh=5.0,Nmr=2.0](n5)(44,-18.0){$u_x$}

    \node[Nw=5.0,Nh=5.0,Nmr=2.0](n6)(64.0,-18.0){$x$}

    \node[Nw=5.0,Nh=5.0,Nmr=2.0](n7)(48.0,-30.0){$y$}

    \drawedge[curvedepth=3,ELside=r](n4,n5){$a$}

    \drawedge[ELside=r](n5,n6){$v$}

    \drawedge[ELside=r](n7,n6){$v$}

    \put(52,-38){in $\Gamma'$}

    \node[Nw=5.0,Nh=5.0,Nmr=2.0](n8)(80.0,-2.0){$z'$}

    \node[Nw=5.0,Nh=5.0,Nmr=2.0](n9)(112.0,-18.0){$x$}

    \node[Nw=8.0,Nh=5.0,Nmr=2.0](n10)(92.0,-18.0){$\scriptstyle y = u_x$}

    \drawedge[curvedepth=3,ELside=r](n8,n10){$a$}

    \drawedge[ELside=r](n10,n9){$v$}
    
    \put(95,-38){in $\Gamma''$}

\end{picture}
\caption{$x\in V_3$, $a\in Y \setminus \{v\}$}
\label{fig xV3}
\end{figure}

Next, we proceed to trimming $\Gamma''$, since we want to compute
$\phi(\Gamma)$, which is equal to $\cc(\Gamma'')$.  The analysis above
shows that $\Gamma''$ has two kinds of endpoints:
\begin{itemize}
    \item vertices $x\in V_2$ such that $\hl_\Gamma(x) \subseteq Y$
    (i.e., $\hl_\Gamma(x) \subseteq Y\setminus \{ v\}$),
    
    \item vertices $x\in V_3$ such that $\hl_\Gamma(x) \subseteq Y$.
\end{itemize}
Suppose first that $x\in V_2$ and $\hl_\Gamma(x) \subseteq Y\setminus
\{ v\}$.  Since $x$ is not an endpoint in $\Gamma$, $\hl_\Gamma(x)$
contains at least 2 elements $a\ne b$ and hence, $\hl_{\Gamma''}(u_x)
= \hl_{\Gamma'}(u_x)$ has at least three elements, see Figure~\ref{fig
xV2trim}.  Removing the vertex $x$ and the only adjacent edges
(labeled $v$ from $u_x$ to $x$, and $\bar v$ from $x$ to $u_x$) leaves
a non-singleton link at $u_x$, that is, trimming $x$ does not create a
new endpoint.

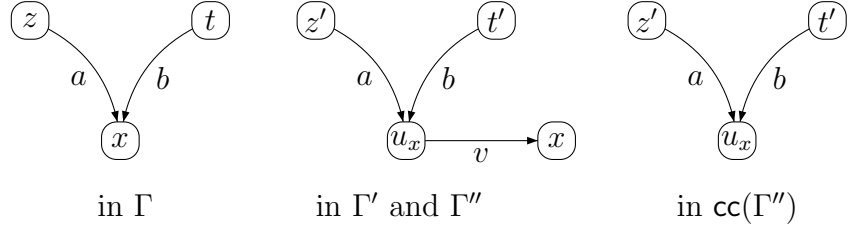
\begin{figure}[ht]
    \centering
\begin{picture}(106,28)(0,-28)

    \node[Nw=5.0,Nh=5.0,Nmr=2.0](n1)(0.0,-2.0){$z$}

    \node[Nw=5.0,Nh=5.0,Nmr=2.0](n2)(24.0,-2.0){$t$}

    \node[Nw=5.0,Nh=5.0,Nmr=2.0](n3)(12.0,-18.0){$x$}

    \drawedge[curvedepth=3.0,ELside=r](n1,n3){$a$}

    \drawedge[curvedepth=-3.0](n2,n3){$b$}
    
    \put(9,-28){in $\Gamma$}

    \node[Nw=5.0,Nh=5.0,Nmr=2.0](n4)(38.0,-2.0){$z'$}

    \node[Nw=5.0,Nh=5.0,Nmr=2.0](n5)(50.0,-18.0){$u_x$}

    \node[Nw=5.0,Nh=5.0,Nmr=2.0](n6)(70.0,-18.0){$x$}

    \node[Nw=5.0,Nh=5.0,Nmr=2.0](n7)(62.0,-2.0){$t'$}

    \drawedge[curvedepth=3.0,ELside=r](n4,n5){$a$}

    \drawedge[ELside=r](n5,n6){$v$}

    \drawedge[curvedepth=-3.0](n7,n5){$b$}

    \put(38,-28){in $\Gamma'$ and $\Gamma''$}

    \node[Nw=5.0,Nh=5.0,Nmr=2.0](n8)(82.0,-2.0){$z'$}

    \node[Nw=5.0,Nh=5.0,Nmr=2.0](n9)(94.0,-18.0){$u_x$}

    \node[Nw=5.0,Nh=5.0,Nmr=2.0](n10)(106.0,-2.0){$t'$}

    \drawedge[curvedepth=3.0,ELside=r](n8,n9){$a$}

    \drawedge[curvedepth=-3.0](n10,n9){$b$}
    
    \put(86,-28){in $\cc(\Gamma'')$}

\end{picture}
\caption{$x\in V_2$, $a,b\in Y \setminus \{v\}$, $a\ne b$}
\label{fig xV2trim}
\end{figure}

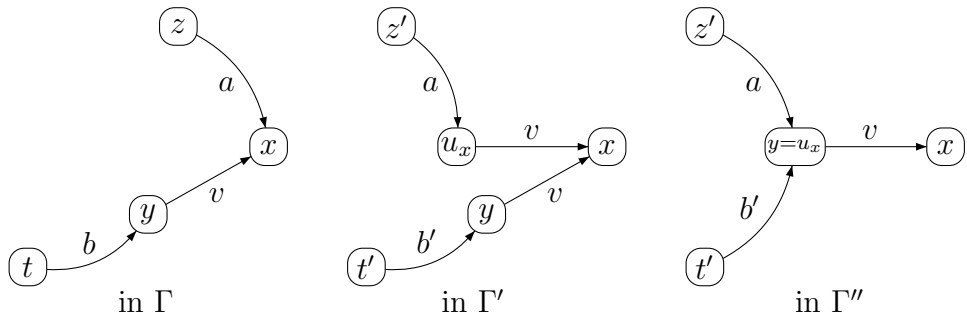
\begin{figure}[ht]
    \centering
\begin{picture}(122,40)(0,-40)

    \node[Nw=5.0,Nh=5.0,Nmr=2.0](n1)(20.0,-2.0){$z$}

    \node[Nw=5.0,Nh=5.0,Nmr=2.0](n2)(16.0,-27.0){$y$}

    \node[Nw=5.0,Nh=5.0,Nmr=2.0](n3)(32.0,-18.0){$x$}

    \node[Nw=5.0,Nh=5.0,Nmr=2.0](n11)(0.0,-34.0){$t$}

    \drawedge[curvedepth=3,ELside=r](n1,n3){$a$}

    \drawedge[ELside=r](n2,n3){$v$}
    
    \drawedge[curvedepth=-3](n11,n2){$b$}

    \put(12,-40){in $\Gamma$}

    \node[Nw=5.0,Nh=5.0,Nmr=2.0](n4)(49.0,-2.0){$z'$}

    \node[Nw=5.0,Nh=5.0,Nmr=2.0](n5)(57,-18.0){$u_x$}

    \node[Nw=5.0,Nh=5.0,Nmr=2.0](n6)(77.0,-18.0){$x$}

    \node[Nw=5.0,Nh=5.0,Nmr=2.0](n66)(61.0,-27.0){$y$}

    \node[Nw=5.0,Nh=5.0,Nmr=2.0](n7)(45,-34){$t'$}

    \drawedge[curvedepth=3,ELside=r](n4,n5){$a$}

    \drawedge(n5,n6){$v$}

    \drawedge[ELside=r](n66,n6){$v$}

    \drawedge[curvedepth=-3](n7,n66){$b'$}

    \put(55,-40){in $\Gamma'$}

    \node[Nw=5.0,Nh=5.0,Nmr=2.0](n8)(90.0,-2.0){$z'$}

    \node[Nw=8.0,Nh=5.0,Nmr=2.0](n9)(102.0,-18.0){$\scriptstyle y = u_x$}

    \node[Nw=5.0,Nh=5.0,Nmr=2.0](n99)(122.0,-18.0){$x$}

    \node[Nw=5.0,Nh=5.0,Nmr=2.0](n10)(90.0,-34.0){$t'$}

    \drawedge[curvedepth=3,ELside=r](n8,n9){$a$}

    \drawedge(n9,n99){$v$}

    \drawedge[curvedepth=-3](n10,n9){$b'$}
    
    \put(102,-40){in $\Gamma''$}

\end{picture}
\caption{$x\in V_3$, $a\in Y \setminus \{v\}$,  $b\ne \bar v$, 
$b'\in\{b,v\}$, $b' \in Y^c \cup\{v\}\setminus \{ \bar v\}$}
\label{fig xV3trim}
\end{figure}

Suppose now that $x\in V_3$ and $\hl_\Gamma(x) \subseteq Y$, and let
$y$ be the initial vertex of the $v$-labeled edge of $\Gamma$ into
$x$.  Since $x \in V_3$, there is an $a$-labeled edge of $\Gamma$ into
$x$ (say, from vertex $z\in V$) for some $a\in Y\setminus\{v\}$, and
therefore, $a \in \hl_{\Gamma'}(u_x) \cap (Y \setminus \{v\})$.  Since
$y$ is not an endpoint in $\Gamma$, $\link_\Gamma(y)$ contains an edge
labeled $b \ne \bar v$, see Figure~\ref{fig xV3trim}.  Then
$\link_{\Gamma'}(y)$ contains an edge labeled $b'$ with $b' = v$ if
$b\in Y$, and $b' = b$ otherwise.  In particular, $b' \in Y^c
\cup\{v\}\setminus \{ \bar v\}$, and hence $b'\ne a$.  Finally, the vertices
$u_x$ and $y$ are identified in $\Gamma''$, so $\link_{\Gamma''}(y)$
contains the distinct elements $a,b'$ and again, trimming $x$ does not
create a new endpoint.

It follows that trimming these two families of endpoints suffices to 
yield $\cc(\Gamma'') =\varphi (\Gamma)$.  The total number of
vertices trimmed in this process is $\card(\{x\in V_2 \cup
V_3 \mid \hl_\Gamma(x)\subseteq Y\}) =\card(\{x\in V \mid
\hl_\Gamma(x)\subseteq Y\})$, so
\begin{eqnarray*}
    |\phi(\Gamma )|-|\Gamma | \enspace = \enspace |\cc(\Gamma''
    )|-|\Gamma | &=& |\Gamma'' |-\card(\{x\in V \mid
    \hl_\Gamma(x)\subseteq Y\}) -|\Gamma | \\    
    &=&\card(V_2 )-\card(\{x\in V\mid \hl_\Gamma(x)\subseteq Y\}).
\end{eqnarray*}
In this count, we observe that each vertex $x\in V$ may contribute
positively (if $x\in V_2$) and negatively (if $\hl_\Gamma(x)\subseteq
Y$), that is
$$|\phi(\Gamma)| - |\Gamma| = \sum_{x\in V} \delta(x),$$
where
$$\delta(x) = \cases{
      +1 & if $x\in V_2$ and $\hl_\Gamma(x) \not\subseteq Y$,\cr
      -1 & if $x\not\in V_2$ and $\hl_\Gamma(x) \subseteq Y$,\cr
      0 & otherwise.}
$$
This is equivalent to
$$
\delta(x) = \cases{
      +1 & if $\hl_\Gamma(x)$ meets $Y$ and $Y^c$, and
      $v\not\in\hl_\Gamma(x)$,\cr
      -1 & if $\hl_\Gamma(x) \subseteq Y$ and $v\in\hl_\Gamma(x)$,\cr
      0 & otherwise,}
$$
which yields the expected formula, $|\phi(\Gamma)| - |\Gamma| =
\capac_{W_\Gamma}(Y) - \deg_{W_\Gamma}(v)$.
\eop

\begin{figure}[ht]
	\centering
    \begin{picture}(110,48)(0,-48)

	\node[Nw=5.0,Nh=5.0,Nmr=2.0](n1)(2.0,-0.0){1}

	\node[Nw=5.0,Nh=5.0,Nmr=2.0](n2)(2.0, -20.0){2}

	\node[Nw=5.0,Nh=5.0,Nmr=2.0](n3)(2.0,-40.0){3}

	\node[Nw=5.0,Nh=5.0,Nmr=2.0](n4)(25.0,-0.0){4}

	\node[Nw=5.0,Nh=5.0,Nmr=2.0](n5)(25.0,-20.0){5}

	\node[Nw=5.0,Nh=5.0,Nmr=2.0](n6)(25.0,-40.0){6}

	\drawedge(n1,n4){$a$}

	\drawedge[ELside=r](n1,n2){$b$}

	\drawedge[ELside=r](n2,n3){$d$}

	\drawedge(n5,n3){$c$}

	\drawedge(n5,n6){$d$}

	\drawedge(n3,n6){$e$}

	\drawedge(n4,n5){$a$}

	\drawedge[curvedepth=3.0](n2,n5){$c$}

	\drawedge[curvedepth=3.0](n5,n2){$a$}
	
	\put(10,-48){$\Gamma$}

	\node[Nw=5.0,Nh=5.0,Nmr=2.0](n11)(42.0,-0.0){1}

% 	\node[Nw=5.0,Nh=5.0,Nmr=2.0](n12)(42.0, -20.0){2}

	\node[Nw=5.0,Nh=5.0,Nmr=2.0](n13)(42.0,-40.0){3}

	\node[Nw=5.0,Nh=5.0,Nmr=2.0](n14)(65.0,-0.0){4}

	\node[Nw=5.0,Nh=5.0,Nmr=2.0](n15)(65.0,-20.0){5}

	\node[Nw=5.0,Nh=5.0,Nmr=2.0](n16)(65.0,-40.0){6}
	
	\drawedge(n11,n14){$a$}

	\drawedge(n11,n15){$b$}

	\drawedge[ELside=r](n14,n15){$a$}

	\drawedge(n15,n13){$d$}

	\drawedge(n13,n16){$e$}

	\drawedge[curvedepth=6.0](n14,n16){$d$}

	\drawedge[ELside=r](n14,n13){$c$}

	\drawloop[loopangle=-90.0,loopdiam=6](n15){$c$}

	\put(50,-48){$\phi(\Gamma)$}

	\node[Nw=5.0,Nh=5.0,Nmr=2.0](n21)(87.0,-0.0){1}

	\node[Nw=5.0,Nh=5.0,Nmr=2.0](n27)(82.0,-10.0){7}

	\node[Nw=5.0,Nh=5.0,Nmr=2.0](n22)(87.0, -20.0){2}

	\node[Nw=5.0,Nh=5.0,Nmr=2.0](n23)(87.0,-40.0){3}

	\node[Nw=5.0,Nh=5.0,Nmr=2.0](n24)(110.0,-0.0){4}

	\node[Nw=5.0,Nh=5.0,Nmr=2.0](n25)(110.0,-20.0){5}

	\node[Nw=5.0,Nh=5.0,Nmr=2.0](n26)(110.0,-40.0){6}
	
	\drawedge(n21,n24){$a$}

	\drawedge(n24,n25){$a$}

	\drawedge[ELside=r](n25,n22){$a$}

	\drawedge[ELside=r,curvedepth=6.0,ELpos=65](n24,n23){$c$}

	\drawedge(n25,n26){$d$}

	\drawedge[ELside=r](n22,n23){$d$}

	\drawedge(n23,n26){$e$}

	\drawedge(n27,n21){$a$}

	\drawedge[ELside=r](n27,n22){$b$}

	\drawloop[loopangle=-45.0,loopdiam=6](n25){$c$}
	
	\put(95,-48){$\psi(\Gamma)$}

    \end{picture}
    \caption{The graphs $\Gamma$, $\phi(\Gamma)$ and $\psi(\Gamma)$}
    \label{ex lemma}
    \end{figure}
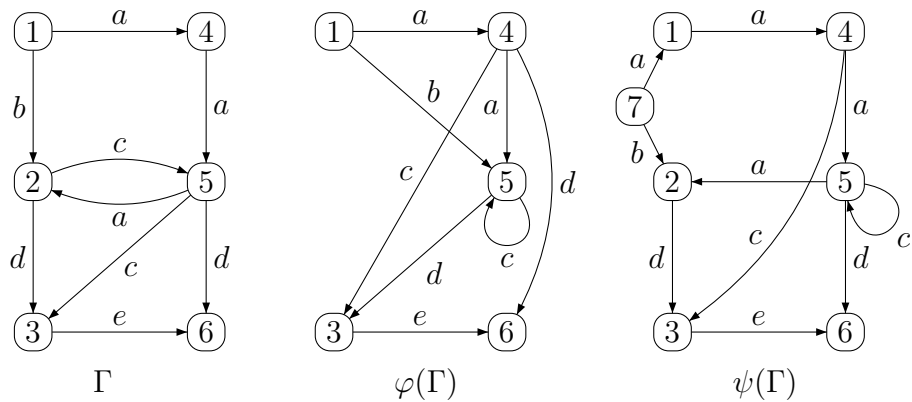

\begin{example}
    Consider the 6-vertex $A$-graph $\Gamma$ in Figure~\ref{ex lemma}.
    Its Whitehead hypergraph $W_\Gamma$ was computed in
    Example~\ref{ex Whitehead graph} and we have $\deg_{W_\Gamma}(a) =
    3$.  Let $\phi$ be the Whitehead automorphism specified by the
    pair $(a,\{a, b, \bar c, \bar d\})$.  We note that
    $\capac_{W_\Gamma}(\{a, b, \bar c, \bar d\}) = 2$, and that
    $\phi(\Gamma)$ then has $6 + 2 - 3 = 5$ vertices, in conformity
    with Proposition~\ref{prop technical}.  The graph $\phi(\Gamma)$
    is shown in Figure~\ref{ex lemma}, as well as the graph
    $\psi(\Gamma)$, where $\psi$ is the Whitehead automorphism
    specified by the pair $(a, \{a, \bar b, \bar c\})$.  Since
    $\capac_{W_\Gamma}(\{a, \bar b, \bar c\}) = 4$ , $\psi(\Gamma)$
    must have size $6 + 4 - 3 = 7$ by Proposition~\ref{prop
    technical}.
\end{example}

\begin{remark}\label{rk complement}
    Let $\Gamma$ be a cyclically reduced $A$-graph, let $v\in \tilde
    A$ and let $Y \subseteq \tilde A$.  It is easily verified that $Y$
    is a $v$-cut if and only if its complement $Y^c$ is a $\bar
    v$-cut, $\capac_{W_\Gamma}(Y) = \capac_{W_\Gamma}(Y^c)$ and
    $\deg_{W_\Gamma}(v) = \deg_{W_\Gamma}(\bar v)$ (equal to the
    number of $v$-labeled edges in $\Gamma$).
%     In particular, if $\phi$ and $\psi$ are the Whitehead
%     automorphisms specified respectively by $(v,Y)$ and $(\bar
%     v,Y^c)$, then $|\phi(\Gamma)| = |\psi(\Gamma)|$.  In fact, it is
%     clear that $\phi = \chi_v \circ \psi$.
\end{remark}

\subsection{Relative complexity of the Whitehead minimization problem}\label{sec relative complexity}

The algorithms to solve the WMP discussed in Section~\ref{sec W
theorem}, can now be modified as follows.

We first consider the case of conjugacy classes of finitely generated
subgroups.  Let $H$ be a cyclically reduced subgroup.  First let
$\vec\phi = (\id)$ and $\Gamma = \Gamma(H)$.  Then repeatedly apply
the following steps: compute the Whitehead hypergraph $W_{\Gamma}$ and
for each $v\in A$, find a $v$-cut $Y_v$ of $\tilde A$ that minimizes
$\capac_{W_{\Gamma}}(Y_v)$; if $\min_{v\in
A}(\capac_{W_{\Gamma}}(Y_v)-\deg_{W_{\Gamma}}(v)) < 0$, let $\psi$ be
the Whitehead automorphism specified by $(v,Y_v)$ where $v$ realizes
the above minimum, and replace $\Gamma$ by $\psi(\Gamma)$ and
$\vec\phi$ by $(\psi, \vec\phi)$; otherwise, stop and output $\Gamma$
and $\vec\phi$.  Finally, choose arbitrarily a vertex 1 in $\Gamma$
and use the procedure discussed in Section~\ref{sec subgroups} to
output a basis of $H'$.

The difference with the algorithm in Section~\ref{sec W theorem} lies
in the fact that instead of trying every Whitehead automorphism to
find one that decreases the size of the cyclically reduced $A$-graph,
we directly select one that will yield the maximum size decrease.  The
fact that we consider only $v$-cuts where $v\in A$ is justified by
Remark~\ref{rk complement}.

Passing from the above algorithm to one that solves the WMP for
subgroups, is done as in Section~\ref{sec W theorem}.

In order to estimate the complexity of the reworded algorithm, we let
$g(n,r)$ be the complexity of the following problem, which we call the
\textit{Whitehead hypergraph min-cut problem}:

\smallskip

{\narrower\narrower\noindent if $\card(A) = r$, given $W_\Gamma$, the
Whitehead hypergraph of a cyclically reduced $A$-graph $\Gamma$ of
size $n$ and a letter $v\in A$, find a $v$-cut $Y$ of $\tilde A$
minimizing $\capac_{W_\Gamma}(Y)$.\par}

\begin{fact}\label{fact WMP gnr}
    As we already saw in Fact~\ref{fact WMP conjclass}, the cost of
    the construction of $\Gamma(H)$, if the input is a set of
    generators of $H$ of total length $n$, is $O(n^2\log(nr))$, and
    $\Gamma(H)$ has at most $n$ vertices.  The computation of the
    image of a cyclically reduced graph of size $n$ under a Whitehead
    automorphism takes time $O(n^2r^2\log(nr))$.  Moreover, the
    Whitehead hypergraph of a size $n$ cyclically reduced $A$-graph
    has $n$ hyperedges, and is computed in time $O(nr\log r)$
    (Fact~\ref{fact complexity Wh graph}).  Finding the degree of a
    vertex in such a hypergraph takes time $O(nr\log r)$.
    
    Then the complexity of each iterating step of our algorithm is at
    most $O(nr\log r + r(g(n,r) + nr\log r) + n^2r^2\log(nr)) =
    O(n^2r^2\log(nr) + rg(n,r))$.  Since there are at most $n$
    iterating steps, the complexity of the full algorithm is
    $O(n^2\log(nr) + n(n^2r^2\log(nr) + rg(n,r)) + n^2\log(nr))$, that
    is, $O(n^3r^2\log(nr) + nrg(n,r))$.
\end{fact}

\begin{fact}\label{fact WMP cyclic}
    Let the \textit{Whitehead graph min-cut problem} be an instance of
    the Whitehead hypergraph min-cut problem where the input is the
    Whitehead graph $W_u$ of a cyclic word $u$, and let $g'(n,r)$ be the
    complexity function of this problem.  Reasoning as in
    Fact~\ref{fact WMP gnr}, we find that the complexity of our
    algorithm to solve the WMP for cyclic words is $O(n^2\log r +
    nrg'(n,r))$.
\end{fact}

%%%%%%%%%%%%%%%%%%%%%%
\section{Main result}\label{sec main}

To conclude our work, we need to find algorithms to solve the
Whitehead hypergraph min-cut problem and its graph analogue in time
polynomial in $n$ and $r$.

%%%%%%%%%%%%%%%%%%%%%%
\subsection{On minimizing the capacity of a cut}\label{sec mincut}

%%%%%%%%%%%%%%%%%%%%%%
\subsubsection{The general case}

The solution of the Whitehead hypergraph min-cut problem can be
reduced to a standard problem in combinatorial optimization, that of
the minimization of submodular functions.  A real-valued function $f$,
defined on the powerset of a set $B$, is said to be
\textit{submodular} if $f(Y\cup Z) + f(Y \cap Z) \le f(Y) + f(Z)$ for
any $Y,Z\subseteq B$.  We first verify the following fact.

\begin{lemma}\label{lemma submodular}
    Let $W = (B,D,\kappa)$ be a hypergraph.  The map $Y \longmapsto
    \capac_W(Y)$, defined on the powerset of $B$, is submodular.
\end{lemma}

\proof
Let $d$ be a hyperedge of $W$.  The contribution of $d$ to
$\capac_W(Y)$ is $0$ if $\kappa(d)\subseteq Y$ or $\kappa(d) \subseteq
Y^c$, and 1 otherwise.  In particular:
\begin{itemize}
    \item if the contribution of $d$ to $\capac_W(Y\cup Z) +
    \capac_W(Y \cap Z)$ is 2, then $\kappa(d)$ meets $Y\cup Z$,
    $(Y\cup Z)^c$, $Y\cap Z$ and $(Y\cap Z)^c$; then it meets
    $Y$, $Y^c$, $Z$ and $Z^c$, so that the contribution of $d$ to
    $\capac_W(Y) + \capac_W(Z)$ is 2 as well;
    
    \item if the contribution of $d$ to $\capac_W(Y) + \capac_W(Z)$ is
    0, then $\kappa(d)$ is contained in $Y$ or in $Y^c$, and it is
    contained in $Z$ or in $Z^c$; equivalently, it is contained in
    $Y\cap Z$, $Y^c\cap Z$, $Y\cap Z^c$ or $Y^c\cap Z^c$; in
    particular, the contribution of $d$ to $\capac_W(Y\cup Z) +
    \capac_W(Y \cap Z)$ is 0 as well;
    
    \item in all other cases, the contribution of $d$ to
    $\capac_W(Y\cup Z) + \capac_W(Y \cap Z)$ is at most 1 and its
    contribution to $\capac_W(Y) + \capac_W(Z)$ is at least 1.
\end{itemize}
This concludes the proof.
\eop

Let $v$ and $W_\Gamma$ be an instance of the Whitehead hypergraph
min-cut problem.  For each $Y\subseteq V\setminus\{v,\bar v\}$, let
$f(Y) = \capac_{W_\Gamma}(Y\cup\{v\})$.  We note that $Y$ minimizes
function $f$ if and only if $Y\cup \{v\}$ has minimum capacity among
the $v$-cuts of $\tilde A$.  It is easily derived from
Lemma~\ref{lemma submodular} that $f$ is submodular.

It follows from results of Gr\"otschel, Lov\'asz and Schrijver
\cite{GLS} that $f$ can be minimized by an algorithm that makes a
polynomial (in $r$) number of oracle calls (queries to evaluate $f$ on
a given argument).  In our situation, given a subset $Y\subseteq
V\setminus\{v,\bar v\}$, computing $f(Y) =
\capac_{W_\Gamma}(Y\cup\{v\})$ takes time $O(nr\log r)$, so the Whitehead
hypergraph min-cut problem can be solved in time polynomial in $n$ and
$r$.  According to Queyranne \cite[pp.  3-4]{Queyranne}, the number of
oracle calls is $O(r^4)$, so the running time of the algorithm is
$O(nr^5\log r)$.

A more recent result of Cunningham \cite{Cunningham} gives a
minimization algorithm with running time $O(Mr^3\log(Mr))$, where $M$
is an upper bound on the maximum value of $f$.  In our case, the value
of $f$ is at most the number of hyperedges in $W_\Gamma$, namely $n$,
so Cunningham's algorithm runs in time $O(nr^3\log(nr))$.  We may take
$g(n,r) = nr^3\log(nr)$.

\begin{remark}\label{rk algos hypergraph mincut}
    The efficient minimization of submodular functions is an active
    research topic, and more recent work offers different algorithms
    which can be used for our purpose just as well as Cunningham's.
    We refer the reader for instance to Iwata, Fleischer and Fujishige
    \cite{IWF}, Schrijver \cite{Schrijver} and Iwata~\cite{Iwata}.
\end{remark}

%%%%%%%%%%%%%%%%%%%%%%
\subsubsection{The graph case}

In the case where the cyclically reduced graph $\Gamma$ is a cyclic
word, the Whitehead hypergraph $W_\Gamma$ is in fact a graph.  The
Whitehead graph min-cut problem is a particular case of the more
general \textit{min-cut problem}, also a standard problem in
operational research, for the solution of which there exists a vast
literature.

In its generality, the min-cut problem for graphs is the following.
We are given a directed graph $G = (V,E)$ with vertex set $V$ and edge
set $E$, and a pair $(s,t)$ of distinct vertices of $G$.  In this
problem, there may be several edges from a vertex $x$ to a vertex $y$.
An \textit{$(s,t)$-cut} of $G$ is a subset $Y$ of $V$, containing $s$
and avoiding $t$.  The \textit{capacity} $\capac_G(Y)$ of such a set
is equal to the number of edges that start in $Y$ and end in the
complement of $Y$.  The min-cut problem consists in finding an
$(s,t)$-cut $Y$ that minimizes $\capac_G(Y)$.

There are many algorithms to efficiently solve the min-cut problem,
see below.  In order to solve the Whitehead min-cut problem on
instance $W_u$ and $v$ (see Section~\ref{sec relative complexity}), we
may first turn $W_u$ into a directed graph $W^+_u$ as follows: we
replace each undirected edge between vertices $x$ and $y$ by a pair of
directed edges, one from $x$ to $y$ and the other from $y$ to $x$.
Next, we observe that a $v$-cut in the sense of Section~\ref{sec W
automorphisms} is a $(v,\bar v)$-cut in the sense of the min-cut
problem, and conversely.  Finally, we verify that if $Y$ is a $(v,\bar
v)$-cut, then both notions of capacity of $Y$ coincide, that is,
$\capac_{W_u}(Y) = \capac_{W^+_u}(Y)$.

Thus a $(v,\bar v)$-cut with minimum capacity in $W^+_u$ is also a
$v$-cut with minimum capacity in $W_u$.  In particular, we may take
$g'(n,r)$ to be the time complexity of any algorithm solving the
min-cut problem in a directed graph with $n$ edges and $2r$ vertices.

Finally, we note that Dinic's algorithm solves the min-cut
problem in time $O(nr^2)$ \cite{Dinic}, see \cite[p.  97]{Kozen}, that 
is, we may take $g'(n,r) = nr^2$.

\begin{remark}\label{rk algos mincut}
    There are many polynomial time algorithms to solve the min-cut
    problem, and we refer the reader to Kozen's book \cite[Chaps.
    15--17]{Kozen} for a review of some of those algorithms that rely
    on the max-flow min-cut theorem (Ford and Fulkerson \cite{FF}),
    that is, that consist in maximizing a flow function associated
    with the graph.  Dinic's algorithm mentioned above falls in that
    category.  We note also that Galil's more recent algorithm
    \cite{Galil} works in time $O(n^{2/3}\,r^{5/3})$.
\end{remark}

%%%%%%%%%%%%%%%%%%%%%%
\subsection{Fully polynomial algorithms}

Putting together the results of Sections~\ref{sec relative complexity}
and~\ref{sec mincut}, we get the expected 
result.

\begin{theorem}
    One can solve the WMP in time polynomial in the size $n$ of the
    input and the rank $r$ of the ambient free group.
\end{theorem}    
 
More precisely, on the basis of Facts~\ref{fact WMP gnr} and~\ref{fact
WMP cyclic}, the discussion in Section~\ref{sec mincut} implies the 
following.

\begin{fact}
    The WMP for finitely generated subgroups and for conjugacy classes
    of finitely generated subgroups can be solved in time
    $O((n^2r^4+n^3r^2)\log(nr))$, where $n$ is the size of the input and
    $r = \rank(F)$.
\end{fact}

\begin{fact}
    The WMP for words and for cyclic words can be solved in time
    $O(n^2r^3)$, where $n$ is the size of the input and $r = \rank(F)$.
\end{fact}    

Our main concern in this paper is the fact that the above complexity
functions are polynomial in $n$ and $r$, and we are less concerned
with the exact polynomial that can be achieved.  In fact, we have
phrased our algorithms in a modular way: an algorithm solving the
Whitehead (hypergraph or graph) min-cut problem is called by our
algorithm, and any improvement in the efficiency of the computation of
a min-cut leads to an improvement in the efficiency of our algorithm.
    
It is also worth noting that in the input of the WMP, we may assume
$r\le n$.  Indeed, letters of $A$ that do not occur in the input word
or subgroup may be ignored, for instance by restricting ourselves to
Whitehead automorphisms that fix them (say, leaving them and their
inverses outside any $v$-cut).  This implies immediately the following
more compact results.

\begin{fact}
    The WMP for finitely generated subgroups and for conjugacy classes
    of finitely generated subgroups can be solved in time $O(n^6\log
    n)$, where $n$ is the size of the input, independently of the rank
    of the ambient free group.
\end{fact}    
    
\begin{fact}
    The WMP for words and for cyclic words can be solved in time
    $O(n^5)$, where $n$ is the size of the input, independently of the
    rank of the ambient free group.
\end{fact}    
    
%%%%%%%%%%%%%%%%%%%%%%
\subsection{Consequences}

Recall that a subgroup $H\lefg F(A)$ is a \textit{free factor} of $F$
if any of its bases can be extended to a basis of $F$.  The
\textit{free factor problem} consists in deciding, given $H$, whether
$H$ is a free factor of $F$.  It is immediate that this is the case if
and only if the minimum size of an element of the automorphic orbit of
$H$ is 1.  Therefore we obtain the following corollary.

\begin{corollary}\label{corol ff}
    There is an algorithm that decides the free factor problem (for a
    subgroup given by a set of generators of total length $n$ in a
    rank $r$ free group) in time polynomial in both $n$ and $r$.
\end{corollary} 

It is interesting to compare this result with that obtained by Silva
and Weil \cite{SilvaWeil}.  These authors give a purely
graph-theoretic algorithm to solve the free factor problem on input
$H$ in time $O(n^{2d+2}\log(nr))$, where $d = r - \rank(H)$.
According to the theoretical complexity functions, the result in
Corollary~\ref{corol ff} is stronger in general, but Silva and Weil's
algorithm may be more efficient on large size, large rank inputs.
Computer experiments might be interesting, especially as the latter
algorithm is simpler to implement, and might yield smaller constants.

A word $u\in F(A)$ is \textit{primitive} if it is an
element of some basis of $F(A)$. That is, $u$ is primitive if and 
only if $\langle u\rangle$ is a free factor of $F$. So we also have the 
following corollary.

\begin{corollary}
    There is an algorithm that decides primitivity (of a word of
    length $n$ in a rank $r$ free group) in time polynomial in both
    $n$ and $r$.
\end{corollary} 

Observe that the formula proved in Proposition~\ref{prop technical} is
additive, in the following sense: if $\mathbf{\Gamma} =
(\Gamma_1,\ldots,\Gamma_m)$ is a tuple of cyclically reduced
$A$-graphs, and if $W_{\mathbf{\Gamma}}$ is the Whitehead hypergraph
of this tuple (the union of the $W_{\Gamma_i}$), if $v\in \tilde A$,
$Y$ is a $v$-cut of $\tilde A$ and $\phi$ is the Whitehead
automorphism determined by $(v,Y)$, then
$$\sum_{i=1}^m|\phi(\Gamma_i)| - \sum_{i=1}^m|\Gamma_i| =
\capac_{W_{\mathbf{\Gamma}}}(Y) - \deg_{W_{\mathbf{\Gamma}}}(v).$$
This additivity extends to Gersten's theorem (Theorem~\ref{W theorem}
above) as observed in \cite{Gersten}.  That is, if some automorphism
of $F$ reduces the total size of a tuple of cyclically reduced
$A$-graphs, then some Whitehead automorphism does \cite[Corol.
2]{Gersten}, generalizing Whitehead's result \cite[Prop.
I.4.20]{LS}).  Our argument then also carries over to the complexity
of the \textit{Whitehead minimization problem for tuples of conjugacy
classes of finitely generated subgroups} (to find a tuple of conjugacy
classes with minimum total size, in the automorphic orbit of a given
tuple).

\begin{corollary}\label{corollary tuples}
    There is an algorithm that solves the WMP for tuples of conjugacy
    classes of finitely generated subgroups in time polynomial in both
    $n$ (the sum of the sizes of the given conjugacy classes) and $r$
    (the rank of the ambient free group).
\end{corollary}

\begin{corollary}\label{corollary tuples cyclic words}
    There is an algorithm that solves the WMP for tuples of cyclic
    words in time polynomial in both $n$ (the sum of the sizes of the
    given conjugacy classes) and $r$ (the rank of the ambient free
    group).
\end{corollary}

%%%%%%%%%%%%%%%%%%%%%%
\subsection{A few open questions}

\paragraph{A cut-vertex theorem?}
Connectedness and connected components are defined in hypergraphs as
in graphs: two vertices $b,b'$ of a hypergraph $W = (B, D, \kappa)$
are connected if there exist a sequence of hyperedges
$d_1,\ldots,d_\ell$ such that $b \in \kappa(d_1)$, $\kappa(d_i) \cap
\kappa(d_{i+1}) \ne\emptyset$ for all $1 \le i < \ell$, and $b' \in
\kappa(d_\ell)$.  Let $\Gamma$ be a cyclically reduced $A$-graph, and
say that $v\in \tilde A$ is a \textit{cut-vertex} of $W_\Gamma$ if
removing $v$ and the hyperedges adjacent to it, yields a hypergraph
$W'$ with more connected components than $W_\Gamma$.
    
If $v$ is a cut-vertex, then the connected component of $W$ containing
$v$ splits into at least two non-empty connected components when $v$
is removed; let $Y'$ be one of them, not containing $\bar v$, let $Y =
Y' \cup \{v\}$ and let $\phi\in\W(A)$ be specified by $(v,Y)$.  By
definition of a cut-vertex, the hyperedges connecting $Y$ and $Y^c$
form a proper subset of the hyperedges adjacent to $v$, that is,
$\capac_{W_\Gamma}(Y) < \deg_{W_\Gamma}(v)$.  It follows that
$|\phi(\Gamma)| < |\Gamma|$.  In particular, if $\Gamma = \Gamma(H)$,
then the size of $H$ is not minimal.
    
This generalizes to subgroups a simple part of Whitehead's celebrated
cut-vertex theorem: if $u$ is a cyclically reduced word and $W_u$ has
a cut-vertex, then the length of $u$ is not minimal.  It is known that
the converse does not hold, that is, the Whitehead graphs of some
non-minimal words present no cut-vertex (the word $u=abbaab\in
F(\{a,b\})$ provides an example).  However, the Whitehead cut-vertex
Lemma (see Stallings \cite[Theorem 2.4]{handlebodies}) states the much
deeper result that if $u$ is primitive then its Whitehead graph $W_u$
is either disconnected (in which case some conjugate of $u$ is
contained in a proper free factor of $F$, \cite[Prop.
2.2]{handlebodies}) or it has a cut-vertex.  It would be interesting
to find an analogous statement for subgroups (and cut-vertices in the
Whitehead hypergraph), and to see whether these statements can be
derived from the combinatorial arguments discussed here.

\paragraph{The hard part of the equivalence problem}
In the so-called hard-part of Whitehead's algorithm to solve the
equivalence problem, say for cyclic words, one considers two cyclic
words $u$ and $v$ of minimum length in their automorphic orbit, and
one needs to decide whether a sequence of Whitehead automorphisms
takes $u$ to $v$ without ever changing the cyclic length.  At first
sight, this might require exploring all words of length $|u|$, which
yields an algorithm that is exponential in both $n$ and $r$.
Myasnikov and Shpilrain \cite{MS} (see also Khan \cite{BK}) establish
that the complexity is in fact polynomial in $n$ for $r = 2$, and more
recent work by Lee \cite{DL} suggests that this is probably true for
all values of $r$.  Kapovich, Schupp and Shpilrain \cite{KSS} develop
a remarkable study of the generic-case complexity of this problem.
One might hope to use our method to get rid of the exponential
dependency in $r$ as well.  This would require being able to find
\textit{all} minimal cuts in the Whitehead graph of a cyclic word
\textit{of minimal length}, and it would be interesting to investigate
whether that can be done in polynomial time.

Extending that investigation to minimal cuts in Whitehead hypergraphs 
would naturally be equally interesting.

\paragraph{Is the greedy algorithm the optimal size reduction 
technique?}
In the Whitehead minimization problem, one may be interested in
minimizing the number of Whitehead automorphisms one needs to apply to
a conjugacy class $[H]$ in order to find a minimum size element of its
orbit.  The algorithm discussed in this paper follows the so-called
\textit{greedy} paradigm: at each step of the iteration, one chooses a
Whitehead automorphism that maximizes the size decrement.  This does
not a priori imply that the number of steps is minimized.  It would be
interesting to verify whether such a greedy algorithm is in fact
optimal also in the number of steps, and to get estimates of that
number of steps.

%%%%%%%%%%%%%%%%%%%%%%%%

%%%%%%%%%%%%%%%%%%%%%%%%
%%%%%%%%%%%%%%%%%%%%%%%%
%%%%%%%%%%%%%%%%%%%%%%%%
\end{document}